\documentclass[preprint,12pt]{elsarticle}
 \pdfoutput=1 
\usepackage[utf8]{inputenc}
\usepackage{graphicx,epstopdf,subfigure,mathtools,mathrsfs} 
\usepackage[normalem]{ulem}
\usepackage{numcompress}

\PassOptionsToPackage{usenames,dvipsnames}{xcolor}
\usepackage[usenames,dvipsnames]{xcolor}

\newcommand{\norm}[1]{\left\lVert #1 \right\rVert} 
\newcommand{\mat}[1]{\mathbf{#1}}
\newcommand{\bm}[1]{\mbox{\boldmath $ #1 $}}    
\renewcommand{\vec}[1]{\bm{#1}}                 

\hypersetup{hidelinks}
\ifpdf
\hypersetup{
  pdftitle={2D force constraints in the method of regularized Stokeslets},
  pdfauthor={Ondrej Maxian and Wanda Strychalski}
}
\journal{Commun. Appl. Math. Comput. Sci.}
\begin{document}
\begin{frontmatter}
\title{2D force constraints in the method of regularized Stokeslets}
\author[OM]{Ondrej Maxian} \ead{om759@cims.nyu.edu}
\address[OM]{Courant Institute of Mathematical Sciences, New York, NY 10012 U.S.A.}
\author[WS]{Wanda Strychalski\corref{cor}}  \ead{wis6@case.edu}
\cortext[cor]{Corresponding author.}
\address[WS]{Department of Mathematics, Applied Mathematics, and Statistics, Case Western Reserve University, Cleveland, OH  44106, U.S.A.}
\begin{abstract}
For many biological systems that involve elastic structures immersed in fluid, small length scales mean that inertial effects are also small, and the fluid obeys the Stokes equations. One way to solve the model equations representing such systems is through the Stokeslet, the fundamental solution to the Stokes equations, and its regularized counterpart, which treats the singularity of the velocity at points where force is applied. In two dimensions, an additional complication arises from Stokes' paradox, whereby the velocity from the Stokeslet is unbounded at infinity when the net hydrodynamic force within the domain is nonzero, invalidating any solutions that use the free space Stokeslet. A straightforward computationally inexpensive method is presented for obtaining valid solutions to the Stokes equations for net nonzero forcing. The approach is based on modifying the boundary conditions of the Stokes equations to impose a mean zero velocity condition on a large curve that surrounds the domain of interest. The corresponding Green's function is derived and used as a fundamental solution in the case of net nonzero forcing. The numerical method is applied to models of cellular motility and blebbing, both of which involve tether forces that are not required to integrate to zero. 
\end{abstract}
\begin{keyword}
fluid-structure interaction \sep Stokes flow \sep Stokes' paradox \sep regularized Stokeslets
\MSC 65M80 \sep 74F10 \sep 92C37
\end{keyword}
\end{frontmatter}
\section{Introduction}

Stokes flow refers to the regime of viscous flow where inertial effects are small, and the Navier-Stokes equations simplify to the Stokes equations. For fluid-structure interaction problems in cell biology, such as an elastic red blood cell membrane deforming in capillary flow \cite{pozrikidis2003numerical,pozrikidis2005axisymmetric}, the small length scales of the cell diameter ($\sim$10 $\mu$m) lead to a small Reynolds number. 
Other important phenomena in cell biology that involve zero Reynolds number flow are cell motility \cite{fauci2006biofluidmechanics,lim2013computational,Vanderlei2011} and microorganism swimming \cite{cisneros2007fluid,lauga2009hydrodynamics}.

Because the Stokes equations are linear, boundary integral and boundary element methods can be used to determine the velocity and pressure fields that come from a collection of forces \cite{pozrikidis1992boundary,pozrikidis2001interfacial}. The velocity field generated from a point force is known as a \textit{Stokeslet}. 
One problem that arises when using the Stokeslet in practice is the singularity at the point where the force arises. For closed interfaces, this singularity is integrable, but careful numerical quadratures are necessary to correctly calculate the velocity and pressure \cite{pozrikidis1992boundary,pozrikidis2001interfacial}. 
In \cite{cortez2001method}, Cortez introduced the method of regularized Stokeslets to overcome the singularities in both the pressure and velocity expressions for forces located at scattered points. Instead of the force being applied at a point, the force is applied over a small ball of radius $\epsilon$. The regularized Stokeslet and pressure expressions are then obtained analytically from the particular function used to represent the small ball. The method of regularized Stokeslets can also be used for closed surfaces, bypassing the associated issues with numerical quadrature \cite{pozrikidis1992boundary}.

It is convenient to model and simulate fluid-structure interaction problems in two dimensional domains where model parameter studies can be conducted in a computationally inexpensive manner. Data visualization is also easier in 2D than in 3D.
The free space Stokes equations in 2D are actually ill-posed because the velocity obtained from the free space Stokeslet is \textit{unbounded} at infinity when there is a nonzero net hydrodynamic force acting within the domain of flow. This contradicts the assumption in the derivation of the Stokeslet that the velocity is zero at infinity and renders the problem ill-posed \cite{cortez2001method,morra2018insights,pozrikidis1992boundary}. 
Numerical simulations of such systems can therefore lead to unphysical spurious velocities. The phenomenon of unbounded velocities in systems with nonzero net force, especially as $||\vec{x}|| \rightarrow \infty$, is usually referred to as Stokes' paradox \cite{Vanderlei2011}. We emphasize that this is a unique feature of Stokes flow in 2D. In 3D, the problem is well-posed; the Stokeslet decays to zero at infinity regardless of the net forcing, so the boundary condition is satisfied and the solution is valid for any collection of forces with bounded magnitude. 

One way to ensure that the 2D velocity is valid is to add conditions to the original system of equations. The most straightforward way to do this is to impose additional boundary conditions within the region of interest. 
The method of regularized Stokeslets was employed in \cite{cortez2001method} to simulate the flow due to a cylinder moving at an imposed velocity. In \cite{ainley2008method,blake1971note,cortez2015general}, the method of images was used to add additional Stokeslets outside of the flow domain that enforce a zero boundary condition near a plane wall. In these approaches, there is an additional constraint on the \textit{velocity} that leads to valid solutions near the immersed objects of interest. However, if no boundary conditions within the flow domain are specified by the model of the physical system (e.g. when modeling flexible fibers in Stokes flow \cite{bouzarth2011modeling, cortez2012slender,tornberg2004simulating}), a different approach must be used. 

One such approach is to enforce a constraint on the \textit{force} rather than the velocity, in particular that the net hydrodynamic force over the entire domain be zero. Sometimes, this constraint comes naturally, such as in models of flagellar swimming \cite{wrobel2016enhanced} or fibers immersed in a background flow \cite{bouzarth2011modeling}. However, the only \textit{a priori} requirement of fluid-structure interaction in Stokes flow is that the hydrodynamic force at a point is exactly balanced by the internal and external forces on the immersed structures \cite{nazockdast2017fast}. In fact, there are many systems with zero Reynolds number that contain force imbalances, including any system that contains tether forces or objects tied to boundaries. For example, Cortez's model of a moving cylinder \cite{cortez2001method} had a nonzero net hydrodynamic force within the domain of flow. In this case, one potential solution is to subtract the mean force from the force at each point, which automatically gives a zero-sum total force. Here we show this approach can result in non-physical, displaced equilibrium states. 


Teran and Peskin \cite{teran2009tether} treated the problem of unbalanced forces within the immersed boundary method \cite{peskin1972flow,peskin2002immersed} by adding a unique, constant, velocity throughout the periodic domain to ensure that the net force is zero for all time. In the formulation from \cite{teran2009tether}, an additional constant velocity is permitted because the equations are simulated on a periodic domain, where the solution is unique up to a constant. In this case, it is required that the net force be zero \cite{batchelor2000introduction}. The net zero force requirement in a periodic IB method presents some challenges. For example, tether forces \textit{must} be introduced in the domain in order to simulate body forces (as the authors did in \cite{teran2009tether} when modeling peristaltic pumping). In order for the immersed structures to remain stationary, the tether spring stiffness must be large, which in turn increases the overall stiffness of the numerical scheme and the cost of the IB method formulation as a whole.

We present a method to simulate models in 2D Stokes flow  with net nonzero forcing using the method of regularized Stokeslets. We accomplish this by surrounding the domain by a large circle and constraining the \textit{mean velocity} on the circle to be zero. Given this boundary condition, we derive the corresponding Green's function 
and show that a mean zero velocity at the large circle can be achieved simply by adding a constant velocity to the free space Stokeslet solution throughout a large domain of flow. In this way, we avoid having to solve a linear system on the large circle (as in \cite{Vanderlei2011}). This observation results in an algorithm that is very straightforward to implement. After presenting our method in Sections 2 and 3, we show in Section 4 how it can be applied to 2D models of cells immersed in viscous fluid. In the process, we compare our formulation to both the explicit zero velocity condition on the large circle, e.g. from \cite{Vanderlei2011}, and the force-free formulation obtained from subtracting the mean force at each point.


\section{Mathematical framework}
The steady Stokes equations in two dimensions are
\begin{align}
\label{eq:mombal}
\mu \Delta \vec{u} - \nabla p & = -\vec{f}\\[4 pt]
\label{eq:incomp}
\nabla \cdot \vec{u} & = 0,
\end{align}
where $\mu$ is the fluid viscosity, $p$ is the pressure, $\vec{u}$ is the fluid velocity, and $\vec{f}$ is the hydrodynamic force, exactly equal to the external applied force that comes from fibers or other structures immersed in the fluid \cite{nazockdast2017fast}. We begin by summarizing the method of regularized Stokeslets \cite{cortez2001method} for computing $\vec{u}$ and $p$ from Eqs.\ \eqref{eq:mombal} and \eqref{eq:incomp}. Then we present the modification for addressing Stokes' paradox. 
\subsection{Method of regularized Stokeslets}
In the method of regularized Stokeslets, a force of strength $\vec{f_0}$ is distributed primarily (but not entirely) over a small ball centered on a point $\bm{x_0}$, so that
\begin{equation}
\label{eq:regforce}
\vec{f}(\vec{x}) = \vec{f_0} \phi_{\epsilon}(\vec{x}-\vec{x_0}).
\end{equation}
The Stokes equations can be solved with the force in Eq.\ \eqref{eq:regforce} to derive the resulting velocity and pressure from a given ``blob'' or ``cutoff'' function $\phi_{\epsilon}$. For example, if 
\begin{equation}
\label{eq:blob}
\phi_{\epsilon}(\bm{x}) = \frac{3\epsilon^3}{2\pi(\norm{\vec{x}}^2+\epsilon^2)^{5/2}},
\end{equation}
then
\begin{equation}
\label{eq:regp}
p^{\epsilon}(\bm{x},\bm{x_0}) = \frac{1}{2\pi}(\bm{f_0} \cdot (\bm{x}-\bm{x_0})) \left(\frac{r_0^2+2\epsilon^2+\epsilon \sqrt{r_0^2+\epsilon^2}}{(\sqrt{r_0^2+\epsilon^2}+\epsilon)(r_0^2+\epsilon^2)^{3/2}}\right) 
\end{equation}
and 
\begin{align}
\vec{u}^\epsilon (\vec{x},\vec{x_0}) = &- \frac{\bm{f_0}}{4\pi\mu}\left(\ln{\left(\sqrt{r_0^2+\epsilon^2}+\epsilon\right)} - \frac{\epsilon (\sqrt{r_0^2+\epsilon^2}+2\epsilon)}{(\sqrt{r_0^2+\epsilon^2}+\epsilon)\sqrt{r_0^2+\epsilon^2}} \right)\nonumber \\[4 pt]
& + \frac{1}{4\pi \mu}(\bm{f_0} \cdot (\bm{x}-\bm{x_0})) (\bm{x}-\bm{x_0}) \frac{\sqrt{r_0^2+\epsilon^2}+2\epsilon}{(\sqrt{r_0^2+\epsilon^2}+\epsilon)^2 \sqrt{r_0^2+\epsilon^2}} \label{eq:regveloc}
\end{align}
are the pressure and velocity that result from the force in Eq.\ \eqref{eq:regforce}, where $r_0=\norm{\bm{x}-\bm{x_0}}$. The derivation of these expressions can be found in \cite{cortez2001method}. Notice that for $r_0 \gg \epsilon$, the standard Stokeslet expressions \cite{pozrikidis1992boundary} are recovered,
\begin{equation}
\label{eq:standpres}
p(\bm{x},\bm{x_0}) = \frac{\bm{f_0} \cdot (\bm{x}-\bm{x_0})}{2\pi r_0^2},
\end{equation}
\begin{equation}
\label{eq:standvel}
\vec{u}(\vec{x},\vec{x_0}) = -\frac{\bm{f_0}}{4\pi \mu} \ln({r_0})+(\bm{f_0} \cdot (\bm{x}-\bm{x_0})) \frac{(\bm{x}-\bm{x_0})}{4 \pi \mu r_0^2}.
\end{equation}
The pressure and velocity resulting from a collection of forces $\vec{f_k}$ spread around a collection of points $\vec{x}_k$ is simply a superposition of the results from Eqs.\ \eqref{eq:regp} and \eqref{eq:regveloc}. It is easy to see that if $\displaystyle{\sum_k \vec{f_k} \neq \vec{0}}$, the velocity in Eq.\ \eqref{eq:regveloc} or \eqref{eq:standvel} is unbounded as $\norm{\bm{x}} \rightarrow \infty$, and the boundary conditions $\vec{u} \rightarrow \vec{0}$ are not satisfied as $\norm{\bm{x}} \rightarrow \infty$.
\subsection{Modification for nonzero net force}
\label{sec:deriveuc}
\begin{figure}
\centering     
\includegraphics[width = 0.5\textwidth]{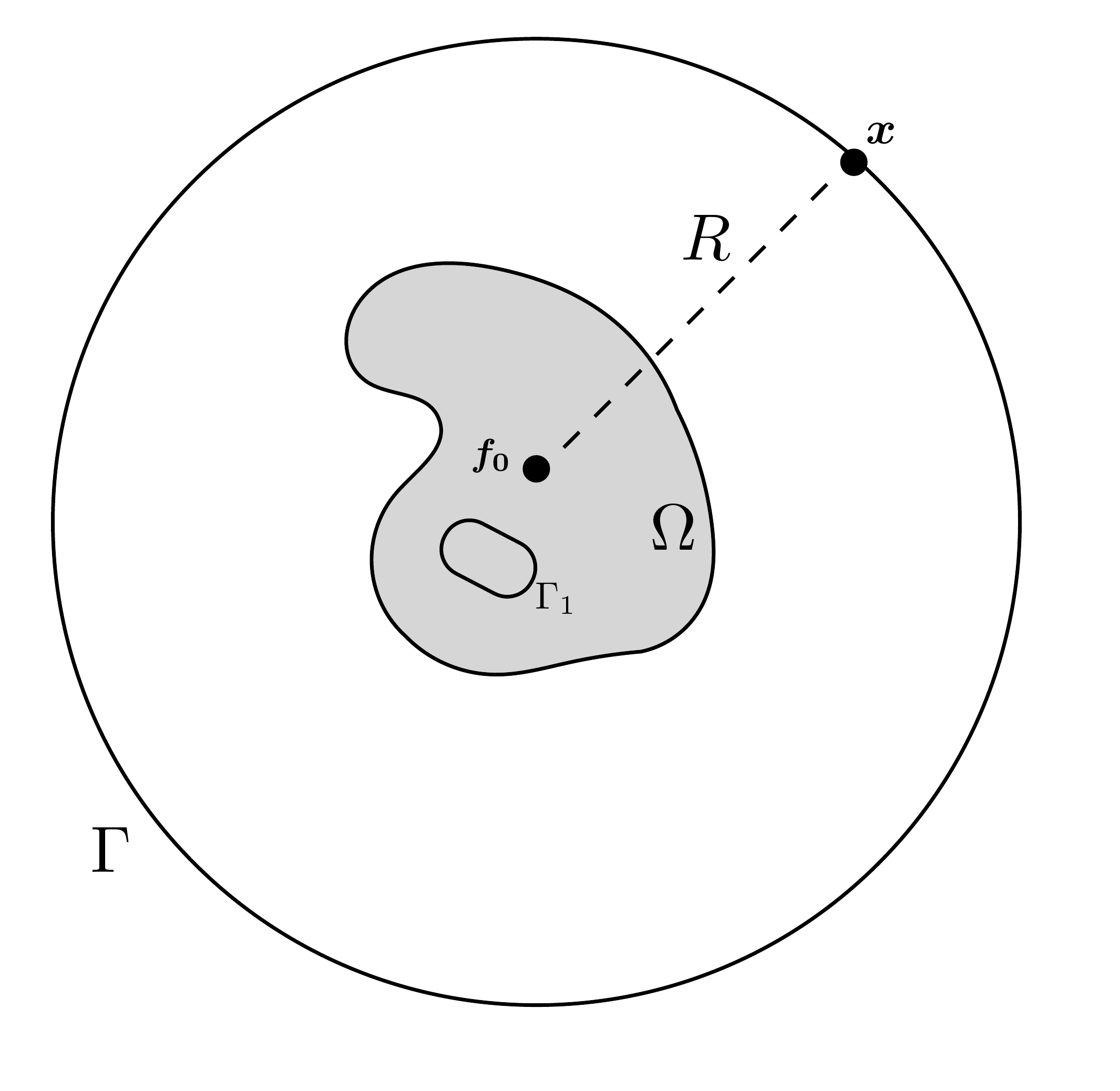}
\caption{$\Omega$ denotes the region bounding immersed interfaces and/or point forces. $\Gamma_1$ indicates an immersed interface and $\vec{f_0}$ denotes a point force at $\vec{x_0}$ enclosed by $\Omega$. A circle of radius $R$ is denoted by $\Gamma$.}
\label{fig:derivation}
\end{figure}

Suppose that all of the forces $\bm{f_k}$ and immersed interfaces in a model system are located within a domain $\Omega$ (see Fig. \ref{fig:derivation}). We note that $\Omega$ is not necessarily an immersed interface, but rather a sort of ``bounding box'' in which all of the forces are contained. If there are no other boundary conditions within $\Omega$, such as a specified velocity on a curve $\Gamma_1$ contained within the domain, the problem is ill-posed and the free space solution for the velocity in Eq.\ \eqref{eq:regveloc} is not valid, even near $\Omega$.
To construct a mathematically valid solution inside some space containing $\Omega$, we surround $\Omega$ with a large circle, denoted by $\Gamma$ with radius $R$ (illustrated in Fig. \ref{fig:derivation}). One approach from \cite{Vanderlei2011} is to enforce a zero velocity boundary condition at every point on the discretized circle. The resulting linear system obtained from Eq.\ \eqref{eq:regveloc} is well conditioned when $\epsilon$ is of magnitude less than or equal to the discrete point spacing on the large circle. The system can be solved for the forces required to obtain a zero velocity on the large circle. These forces can then be used in Eq.\ \eqref{eq:regveloc} to compute the velocity at locations enclosed by the smaller domain $\Omega$. We note that this requires a linear system to be solved at each time value when simulating a model of a dynamic process. 

We take a slightly different approach. Instead of requiring $\vec{u}|_{\Gamma}=\vec{0}$ at every point, we solve Stokes equations with a slightly weaker boundary condition, that the average value of $\vec{u}$ on $\Gamma$, $\langle \vec{u} \rangle |_{\Gamma}=\vec{0}$. Notice that as $R \rightarrow \infty$, the regularized Stokeslet solution in Eq.\ \eqref{eq:regveloc} converges to the free space Stokeslet in Eq.\ \eqref{eq:standvel} and the velocity on the large circle is approximately \textit{constant} because of the dominance of the radially symmetric first term in Eq.\ \eqref{eq:standvel}. If the velocity on the large circle is constant, the two conditions are equivalent. Imposing the mean velocity condition allows us to add an extra velocity $\bm{u^R}$ throughout the domain so that the average velocity on the large circle is zero. 

We begin by deriving the velocity $\bm{u^R}$ in the case of a single point force followed by the generalization to multiple forces by superposition. Let $\bm{f_0}$ be the force at a point $\bm{x_0}$ in $\Omega$. Let $\bm{x}$ be a point on the large circle $\Gamma$ (see Fig.\ \ref{fig:derivation}), and let $r_0 = \norm{\bm{x}-\bm{x_0}}$. Since $\bm{x}$ is on a large circle with arbitrarily large radius $R$, $r_0(\bm{x}) \gg \epsilon$ for all $\bm{x} \in \Gamma$, and we can represent the velocity at the large circle using the standard Stokeslet. Thus the velocity at $\bm{x}$ due to the force $\bm{f_0}$ applied at $\bm{x_0}$ is given by Eq.\ \eqref{eq:standvel}.

Using $s$ as the arclength parameter and treating $r_0=R$ as constant, the average value of $\vec{u}(\vec{x},\vec{x}_0)$ is
\begin{align}
\langle \vec{u}(\vec{x},\vec{x}_0) \rangle & =\frac{1}{2 \pi R} \int_{\Gamma} \left({-\frac{\bm{f_0}}{4 \pi \mu} \ln{(r_0)} + (\bm{f_0} \cdot (\bm{x}-\bm{x_0})) \frac{\bm{x}-\bm{x_0}}{4 \pi \mu r_0^2}} \, \right) \, ds \\[4 pt]
& = \frac{1}{2 \pi R}\left( -\frac{\bm{f_0}}{4 \pi \mu} \ln{R} \int_{\Gamma} ds + \frac{1}{4 \pi \mu R^2}\int_{\Gamma}(\bm{f_0} \cdot (\bm{x}-\bm{x_0})) (\bm{x}-\bm{x_0})\, \, ds \right)\\[4 pt] 
\label{eq:uinter}
& = -\frac{\bm{f_0}}{4 \pi \mu} \ln{R} + \left( \frac{1}{2 \pi R}\right)\frac{1}{4 \pi \mu R^2} \int_{\Gamma} \left(\bm{f_0} \cdot (\bm{x}-\bm{x_0})\right) (\bm{x}-\bm{x_0}) \, \, ds.
\end{align}
The last equality used the fact that $\int_{\Gamma} ds = 2\pi R$. Computing the second integral, we begin by changing to an angle parameterization of $\Gamma$ via $s=R \theta$,
\begin{equation}
\label{eq:arctothet}
\frac{1}{2 \pi R} \int_{\Gamma} \left(\bm{f_0} \cdot (\bm{x}-\bm{x_0})\right) (\bm{x}-\bm{x_0})\, ds = \frac{1}{2 \pi R} \int_0^{2\pi} \left(\bm{f_0} \cdot (\bm{x}-\bm{x_0})\right)(\bm{x}-\bm{x_0}) \, \, R \, d\theta
\end{equation}
Because the circle is rotation invariant, we can place the $x$ axis on the same direction as $\bm{f_0}$ without loss of generality. Therefore, let  $\bm{f_0}=f\begin{pmatrix} 1 \\[2 pt] 0 \end{pmatrix}$. Furthermore, in the limit $R \rightarrow \infty$, $\bm{x}-\bm{x_0}=\bm{x} = \begin{pmatrix} R \cos{\theta} \\[2 pt] R \sin{\theta} \end{pmatrix}$. Then Eq.\ \eqref{eq:arctothet} simplifies to
\begin{align}
\frac{1}{2 \pi R} \int_0^{2\pi} \left(\bm{f_0} \cdot (\bm{x}-\bm{x_0})\right)(\bm{x}-\bm{x_0}) \, \, R \, d\theta & =  \frac{1}{2 \pi} \int_0^{2\pi} fR \cos{\theta} \begin{pmatrix} R \cos{\theta} \\[2 pt] R \sin{\theta} \end{pmatrix} d\theta\\[4 pt]
\label{eq:simplfiedsecond}
& = \frac{R^2 f}{2} \begin{pmatrix} 1 \\[2 pt] 0 \end{pmatrix} = \frac{R^2 \bm{f_0}}{2}. 
\end{align}
Substituting Eq.\ \eqref{eq:simplfiedsecond} into Eq.\ \eqref{eq:uinter}, we have the average velocity over the large circle given the force $\bm{f_0}$ at $\bm{x_0}$ as
\begin{equation}
\label{eq:firstaddedvel}
\langle \vec{u}(\vec{x},\vec{f}_0) \rangle = \frac{\bm{f_0}}{4 \pi \mu} \left(\frac{1}{2} - \ln{R} \right).
\end{equation}
Our goal is to impose a boundary condition on the large circle. Rather than impose a boundary condition pointwise, we impose a weaker condition on the {\em mean velocity} on the large circle, namely that the mean velocity is zero. It follows immediately that this can be done by subtracting the constant velocity in Eq.\ \eqref{eq:firstaddedvel} throughout the domain of flow. The additional velocity due to a point force $\vec{f_0}$ is therefore

\begin{equation}
\label{eq:urptforce}
\vec{u^R}(\bm{f_0})= -\frac{\bm{f_0}}{4 \pi \mu} \left(\frac{1}{2} - \ln{R} \right).
\end{equation}

In the case of multiple time-dependent forces $\bm{f_k}(t)$ (for example, forces that come from an interface such as $\Gamma_1$ in Fig.\ \ref{fig:derivation}), the constant velocity is simply the superposition of velocities from Eq.\ \eqref{eq:urptforce}.
\begin{equation}
\label{eq:addedvel}
\bm{u^R} (t) = \sum_{k=1}^N {\vec{u^R}(\bm{f_k} (t))} = \sum_{k=1}^N {-\frac{\bm{f_k}(t)}{4 \pi \mu} \left(\frac{1}{2} - \ln{R} \right)}.
\end{equation}
This velocity is added throughout the domain of flow to ensure that $\langle \bm{u}(t) \rangle|_{\Gamma} =0$. In effect, the solution from Eq.\ \eqref{eq:regveloc} and Eq.\ \eqref{eq:addedvel} together form the \textit{Green's function for Stokes equations with a mean zero boundary condition on the circle of radius $R$}. We note the Green's function would change if the domain was surrounded by a large square with edgelength $2R$ as opposed to a circle of radius $R$. However, any geometry, such as a square, can be thought of as bounded by two concentric large circles (for a square centered at 0 with edgelength $2R$, the points on the square are between concentric circles of radius $R$ and $R\sqrt{2}$). Because the velocity dependence on $R$ is weak for large $R$ (the derivative of $\vec{u}^{\mathbf{R}}$ scales with $1/R$), altering the geometry of the boundary results in small changes in the velocity on the bounding region.
%

It is no coincidence that Eq.\ \eqref{eq:addedvel} expresses the average value of the velocity due to forces of strength $-\vec{f_k}$. By adding the constant velocity in Eq.\ \eqref{eq:addedvel}, we are effectively adding a force on the large circle that \textit{has the effect of adding an equal and opposite force within the region of interest $\Omega$}. Meanwhile, the addition of a constant velocity throughout the domain results in a relative velocity profile that is unchanged from that computed by Eq.\ \eqref{eq:regveloc}. Our approach contrasts with adding more Stokeslets at arbitrary locations in the domain $\Omega$, which could be problematic because the relative local profile (and subsequent physical conclusions) are dependent on the locations of the additional Stokeslets. However, our approach of imposing forces on the large circle does result in the introduction of a much larger length scale in the problem; the length scale becomes $R$, the radius of the large circle, instead of the length scale of the immersed objects. 


The addition of this constant velocity has no effect on the pressure profile calculated from Eq.\ \eqref{eq:regp}. Because a constant velocity is added, no pressure gradient is generated within the domain. Equivalently, our addition of a constant velocity is a shortcut around explicitly adding forces on the large circle that enforce the zero boundary condition exactly (e.g. \cite{Vanderlei2011}). Due to the nature of the pressure solution in Eq.\ \eqref{eq:regp} (i.e. that it decays as 1/length), 
the additional forces from the boundary condition on the large circle have no effect on the local pressure profile for large $R$. 

We also note that the choice of blob in Eq.\ \eqref{eq:blob} yields the resulting analytical expressions for the regularized pressure and velocity in Eqs. \eqref{eq:regp} and \eqref{eq:regveloc}, respectively.
We present these expressions derived from $\phi_{\epsilon}$ because we use them in our numerical simulations. The derivation of $\bm{u^R}$ is independent of the regularized Stokeslet because it is derived from the true free space Stokeslet. The large circle $\Gamma$ is assumed to be far enough from the domain of interest that the two are equivalent. The method we present here is therefore compatible with any regularization kernel, including compactly supported immersed boundary kernels \cite{bao2016gaussian}. However, IB kernels are generally used over a periodic fluid grid, not over free space, and so the method of \cite{teran2009tether} is more appropriate in that context.

\section{Discretization}
In general, we begin with a collection of $N$ points $\bm{x_k}$ with forces $\bm{f_k}$ in some domain $\Omega$. At each timestep, we compute the velocity of each point $\vec{x_i}$ as
\begin{equation}
\label{eq:totalvel}
\vec{u}(\vec{x_i}) = \bm{u^R}+\sum_{k=1}^N \vec{u}^\epsilon (\vec{x_i},\vec{x_k}),
\end{equation}
where $\bm{u^R}$ is given by Eq.\ \eqref{eq:addedvel} and $\vec{u}^\epsilon (\vec{x_i},\vec{x_k})$ is given by Eq.\ \eqref{eq:regveloc}. Because $\bm{u^R}$ is constant throughout the domain, the calculation in Eq.\ \eqref{eq:totalvel} is $\mathcal{O}(2N+(2N)^2)$ operations. The first $\mathcal{O}(2N)$ operations arise from the computation of the constant additional velocity in Eq.\ \eqref{eq:addedvel}. The second $\mathcal{O}((2N)^2)$ operations come from computing the regularized velocities in Eq.\ \eqref{eq:regveloc} for all of the points.

We compare this operation count to alternative formulations. Suppose that the large circle was discretized with $M$ points and the forces on the large circle solved for explicitly, as in \cite{Vanderlei2011}. This calculation is a $2M \times 2M$ linear solve and requires $\mathcal{O}((2M)^2)$ operations using GMRES or $\mathcal{O}((2M)^3)$ operations if done directly. In addition, the calculation of the added velocities at each point in the domain from the forces on the large circle requires another $\mathcal{O}(2NM)$ flops, and it is unclear how to choose the number and location of the $M$ points. Alternatively, the addition of more Stokeslets in a method similar to the method of images \cite{ainley2008method,blake1971note,cortez2015general} would require $\mathcal{O}(2SN)$ operations to compute the added velocity, where $S$ is the number of added Stokeslets. Our added velocity is computed in $2N$ flops, making it much more efficient than any of these alternatives. 


\subsection{Choosing the radius}
Central to our method is the assumption that the velocity computed from Eq.\ \eqref{eq:standvel} on the large circle $\Gamma$ is relatively constant. The validity of this assumption dictates a lower bound on $R$. In order to test the variation of the velocity on the large circle $\Gamma$, we first impose a force of $\bm{f_0}=\begin{pmatrix} 1\\[1 pt] 0 \end{pmatrix}$ at the origin. Next, we use Eq.\ \eqref{eq:standvel} to compute and measure the velocity from the Stokeslet (i.e. the velocity without the addition of $\bm{u^R}$) on the large circle. 

Specifically, we discretize the circle $\Gamma$ with $N=100$ points and quantify velocity variability by defining
\begin{equation}
\label{eq:velvar}
\sigma_u(R) = \max_i{\left|\frac{u_x^i - \bar{u}_x}{\bar{u}_x}\right|}.
\end{equation}
Here the index $i$ runs from 1 to $N=100$, $u_x^i$ refers to the velocity in the $x$ direction (the direction of the force) at point $i$, and $\bar{u}_x$ refers to the arithmetic mean of $u_x$ taken over $\Gamma$.

\begin{figure}
\centering     
\includegraphics[width = 0.6\textwidth]{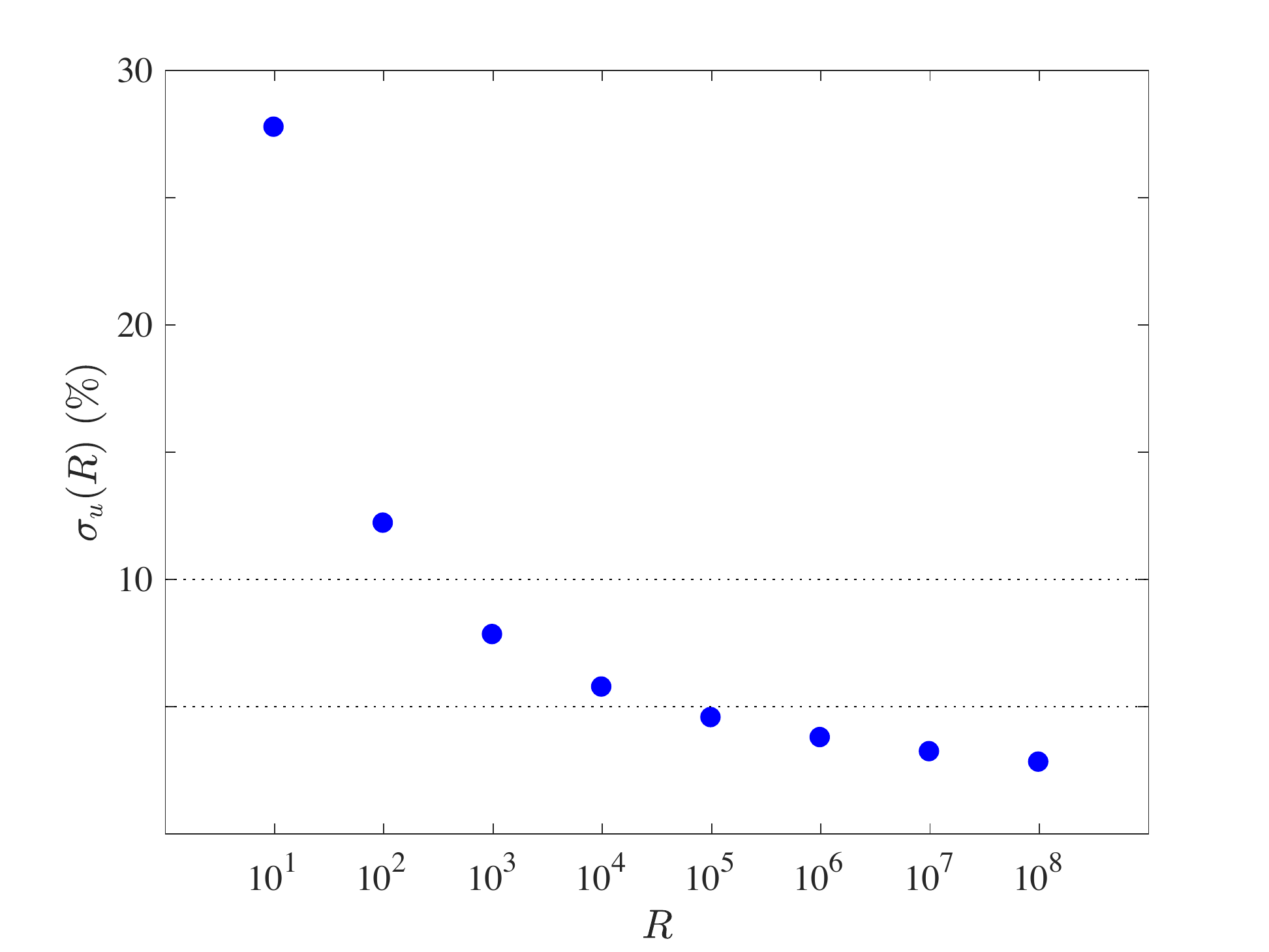}
\caption{Variation of the velocity on the large circle $\Gamma$ for different values of $R$. Dotted lines indicate 10\% and 5\% variation, $\sigma_u(R)$.}
\label{fig:assumptest}
\end{figure}

Fig.\ \ref{fig:assumptest} shows the values of $\sigma_u(R)$ (as a percentage) for different values of $R$. We observe a variation in the velocity $<10\%$ for $R \geq 10^3$ and a variation $<5\%$ for $R \geq 10^5$.

With this in mind, we specify the lower limit on $R$, $R \geq 10^3$. We can find an upper limit on $R$ based on the Reynolds number. Let $v$ and $L$ be the relevant velocity and length scales. Our systems have values of viscosity $\mu$ and density $\rho$ of the same magnitude as water, so that Re$\displaystyle{=\frac{\rho v L}{\mu}=}$ $10^6vL$ is the relevant Reynolds number. Our model systems are from applications in cell biology, so the relevant velocity scale is in $\mu$m/s. Therefore we take $v=10^{-6}$ m/s and Re $=L$, where $L$ is the relevant length scale. In order for Stokes flow to be valid, we need (Re $\ll 1$), which we define to be Re $\leq 0.1$. Then the relevant length scale cannot exceed 0.1 m = $10^5$ $\mu$m. Because we have confined the domain and effectively introduced forces on the large circle, the radius of the large circle is now the largest relevant length scale, and we have determined an upper bound on $R$, $R \leq 10^5$ $\mu$m. We note that this upper bound may change depending on the characteristic time and length scales used to compute the Reynolds number, but it is straightforward to derive it as we have here. Thus we have determined in general that $10^3 \leq R \leq 10^5$. For the examples in Section \ref{sec:exs}, we choose $R=10^3$ $\mu$m to ensure the validity of Stokes equations for these systems. This value of $R$ results in a velocity variation from Eq.\ \eqref{eq:velvar} less than 10\%. In addition, we did not find any additional stiffness when using $R$ is this range (for larger $R$, e.g. $R=10^{16}$, the velocities in Eq.\ \eqref{eq:addedvel} would increase, thereby increasing the overall problem stiffness.


\section{Examples}
\label{sec:unifex}
\label{sec:exs}
The motivation for this work is \textit{tether forces} that arise in the modeling of some biological systems. These are forces that penalize displacement from an initial or resting configuration; points on an immersed object are physically {\em tethered} to other points in 2D space.
We begin by considering a simplified system of tethered particles. This motivating example establishes the need for the additional velocity in Eq.\ \eqref{eq:addedvel}. We then present a model of a cell motility problem where tether forces are useful for modeling the cell's external environment.
We conclude by analyzing a boundary integral model of cellular blebbing with nonzero net forcing that has already been used for modeling bleb initiation and amoeboid cell motility \cite{fang2017combined,lim2012size, lim2013computational}. In all cases, all of the objects are flexible so that no boundary conditions are provided from the physics of each model system.

\subsection{A motivating example}
\label{sec:motex}
\begin{figure}
\label{fig:example1}
\centering     
\subfigure[Initial configuration]{\label{fig:ex1setup}\includegraphics[width=120mm]{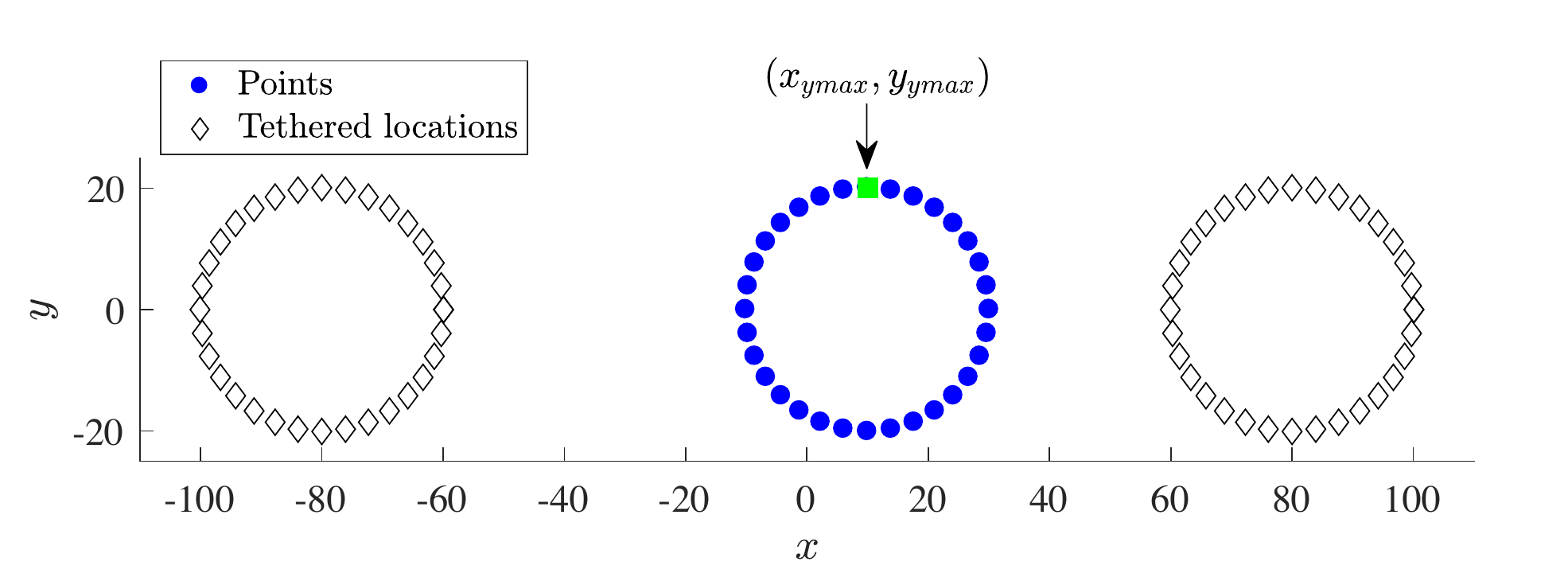}}
\subfigure[Dynamics]{\label{fig:ex1results}\includegraphics[width=80mm]{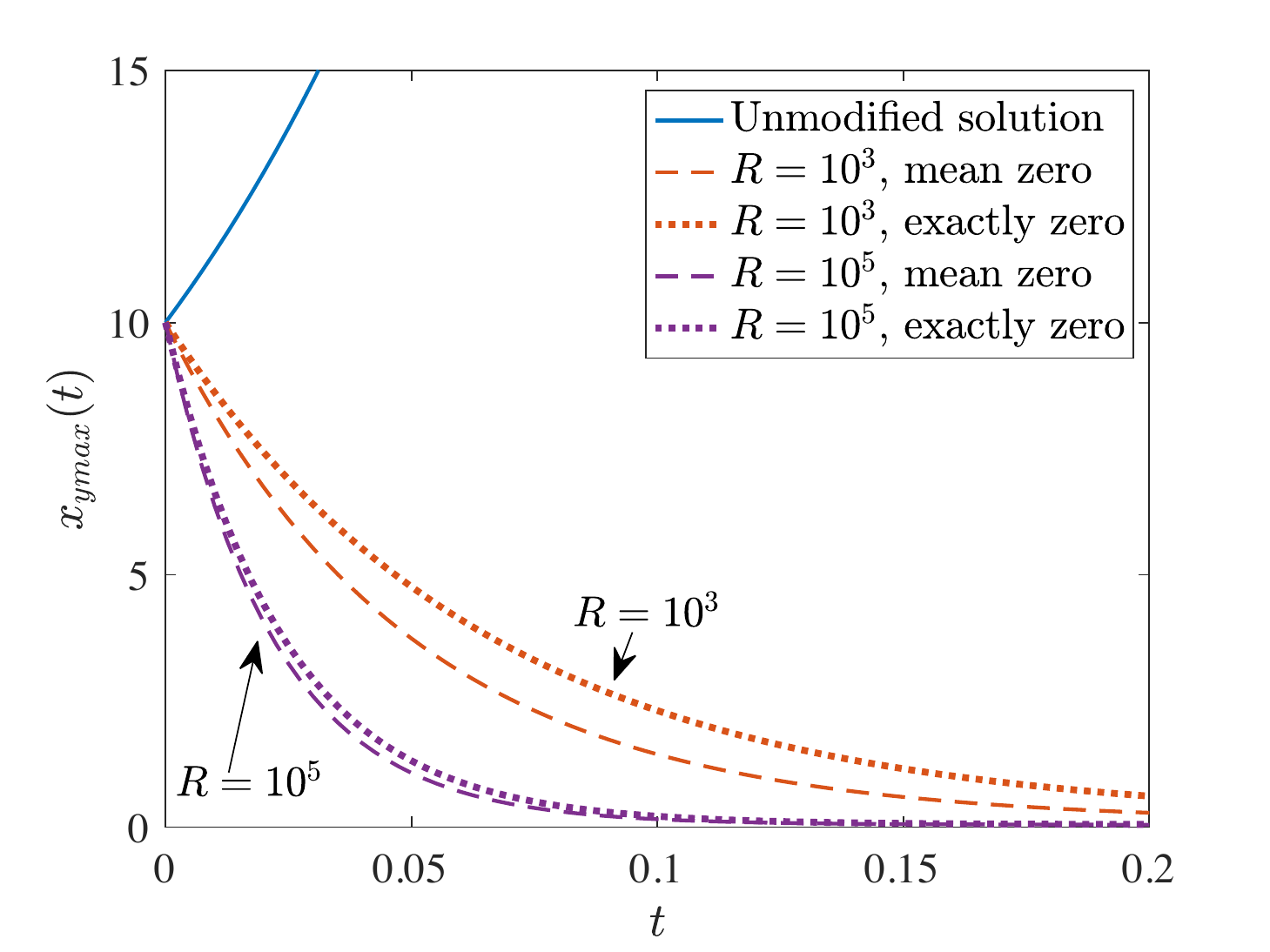}}
\caption{Motivating example of a system of points tethered in place, initially out of equilibrium. (a) Initial configuration of the system. The solid blue dots show the points, while the unfilled black diamonds show the corresponding tethering locations. (b) $x$ coordinate of the point with largest $y$ coordinate, $x_{ymax}$ (filled green square in (a)). Without any corrections, the motion is in the positive horizontal ($+x$) direction independent of number of points and $\epsilon$ (blue line). This nonphysical behavior can be corrected by adding the velocity $\bm{u^R}$ in Eq.\ \eqref{eq:addedvel}. We show data for the values $R=10^3$ (orange) and $R=10^5$ (purple). Dashed lines give the dynamics for the mean velocity subtraction in Eq.\ \eqref{eq:totalvel}. Dotted lines show dynamics when the boundary condition $\vec{u}=\vec{0}$ is exactly enforced on a large circle of radius $R$ discretized with $N=100$ points.}
\end{figure}
We are interested in modeling the motion of a cell through a viscoelastic structure called the extracellular matrix (ECM). For example, our model of the ECM represents a lattice of collagen fibrils immersed in a viscous fluid.
Elasticity of the ECM can be modeled in 2D by a lattice of points that are tied to specified reference points by springs. In order to demonstrate why our methodology is key for modeling this process, we introduce a set of $N=32$ points distributed on a circle of radius $r$, shown in Fig.\ \ref{fig:ex1setup}. The points (solid blue dots) are centered at $(10,0)$ and are tethered to two sets of fixed points (hollow black diamonds) centered at $(80,0)$ and $(-80,0)$. For simplicity, we set
\begin{equation}
\label{eq:tforceex}
    \bm{F}_i (t)=-k_{teth}\left(\bm{X}_i (t)-\bm{X}_i^R + \bm{X}_i(t)- \bm{X}_i^L\right), 
\end{equation}
where $\bm{X}^R_i$ and $\bm{X}^L_i$ stands for the position of point $i$ on the right and left fixed circles, respectively. The parameters for the system are $\mu=1$ Pa-s, $k_{teth}=1$ pN/$\mu$m$^2$, $\epsilon=\frac{2\pi r}{N}=$ the point spacing (although the dynamics are independent of the parameters 
$\mu$, $k_{teth}$, and $\epsilon$). The initial configuration results in a force imbalance, with a net force to the left (negative horizontal direction) on the set of moving points. Physically, we expect the points to move in the negative horizontal direction and approach their equilibrium position exactly between the two fixed fibers, i.e. centered at $x=0$. However, this is not always what occurs when using Eq.\ \eqref{eq:regveloc} to update the velocities. 

Consider Eq.\ \eqref{eq:regveloc}. The dominant part of the first term is \\
$-\ln{r}=$ $-\ln{||\vec{x}-\vec{x_0}||}$, which results in a velocity that goes in the direction \textit{opposite} the force. The second term in Eq.\ \eqref{eq:regveloc} is $\mathcal{O}(1,\epsilon^{-2})$ and contributes to the velocity in the same direction as the force. In order for the dynamics to match our physical intuition, the cumulative contribution of the second term at each point must be greater than that of the first term, so either $\ln{r} \approx 1$ or $\epsilon \ll 1$. Thus, as $r$ becomes large, we expect unphysical behavior. For an $r$ value as small as $r \approx 20$ (determined empirically), no value of $\epsilon$ (larger than machine epsilon) yields physical results. Fig.\ \ref{fig:ex1results} shows the horizontal ($x$) position of the point with the largest $y$ coordinate over time (marked with a green square as $(x_{ymax},y_{ymax})$ in Fig.\ \ref{fig:ex1setup}), where the velocity is computed by Eq.\ \eqref{eq:regveloc}. We observe unphysical motion in the positive horizontal direction (blue line), moving the points to the right and resulting in an increased force imbalance as time increases. While this behavior \textit{dominates} for large values of $r$, it is also present for smaller $r$ and needs to be corrected to give proper, physically correct, simulation results.

Our solution with the additional velocity given in Eq.\ \eqref{eq:addedvel} gives the expected behavior. For the initial configuration in Fig.\ \ref{fig:ex1setup}, the motion of the points computed with the velocity in Eq.\ \eqref{eq:totalvel} shifts the points in the negative horizontal direction (left), allowing them to approach their steady state positions at $x=0$. Specifically, the horizontal component of the velocity at the origin can be decomposed into $\bm{u}^\epsilon \approx 125$ and $\bm{u^R} \approx -325$ when $R=10^3$, where the contribution of $\bm{u^R}$ is necessarily greater to obtain the correct physical motion. Fig.\ \ref{fig:ex1results} shows the horizontal position of the point with the largest $y$ coordinate for values of $R$ that satisfy our derived bounds, $R=10^3$ (orange lines) and $R=10^5$ (purple lines). Dashed lines show the solution obtained from Eq.\ \eqref{eq:totalvel}, and dotted lines show the solution obtained from discretizing the circle of radius $R$ with $N=100$ points and explicitly enforcing $\bm{u}=0$ at those points by determining the additional forces $\bm{f}$ on the large circle via solving a linear system of equations. The linear system of equations in this case has the form $\bm{U}=\mat{M}\bm{F}$, where $\mat{M}$ is a dense $2N \times 2N$ matrix. To get an exact solution, we solve this directly with $LU$ factorization, although it could also be done with GMRES.

We observe that the mean velocity solution gives faster dynamics than the discretized large circle solution. As $R$ increases, the velocity on the large circle approaches a constant value and the solutions approach the same curve (as shown previously in Fig.\ \ref{fig:assumptest}). 
Further, the maximum difference of $1.0 \, \mu$m in the $x$-coordinate between the two curves for $R=10^3$ is only 5.2\% of the smallest system length scale $r$, and for $R=10^5$ this difference decreases to 1.8\%. 

The choice for the value of $R$ in Eq.\ \eqref{eq:addedvel} affects the velocity and the dynamics of the system. However, we note that additional forces solved for by enforcing additional boundary conditions on any geometry $C$, i.e. $\vec{u}|_{C}= \vec{0}$, would also affect dynamics of the system.
For the values of $R$ within the range $10^3 \leq R \leq 10^5$, the \textit{steady-state} behavior of our model and relevant timescales are shown in Fig.\ \ref{fig:ex1results} to be nearly identical, with the timescales differing by about a factor of 2. If accurate transient results are desired, the value of $R$ can be tuned to give dynamics that fit within the relevant timescales, with the caveat that increasing $R$ leads to a larger length scale.

We note another feature of this example: the initial force on the configuration shown in Fig.\ \ref{fig:ex1setup} is \textit{uniform across all of the points}. Suppose one wanted to treat the force imbalance in this example by subtracting the mean force from each point, so the total force sums to zero. Because each point has the same force on it, subtracting the mean force gives zero force and zero velocity at every point. We have therefore shown that subtracting the mean force can create artificial equilibrium configurations. This phenomenon occurs not just in this simple example, but also in a more complicated model of cell motility as discussed in Section \ref{sec:cellmot}. Thus, while subtracting the mean of the forces maintains the relevant system length scales, 
doing so can introduce errors in the resting position of the system. Which avenue to choose in this trade-off is application dependent. 

\subsection{Model of cell motility}
\label{sec:cellmot}
Cell motility is an essential process for wound healing, cancer metastasis, and immune responses \cite{rubinstein2005multiscale}. In three dimensions, a cell can utilize multiple mechanisms to migrate through the surrounding extracellular matrix (ECM)  \cite{tozluoglu2015cost, tozluouglu2013matrix, zhu2016comparison}. The ECM is a dense network of collagen fibers (see \cite{even2005cell}, Fig. 1). Previous studies have used 2D agent based/finite element models \cite{tozluoglu2015cost,tozluouglu2013matrix} to study the effectiveness of bleb-based and protrusion-based mechanisms in different ECM environments. In \cite{zhu2016comparison}, the authors simulated a variety of mechanisms on a cell with a rigid nucleus via force balance equations. Our goal is to extend this model to a flexible nucleus, where the fluid-structure coupling is treated explicitly via the method of regularized Stokeslets. 

We focus here on the movement of the cell via a finger-like protrusive mechanism \cite{legant2010measurement,zhu2016comparison}. In this mechanism, the elastic cell cortex generates random actin-based protrusions. The cortex is the thin layer of the actin cytoskeleton that is attached to the cell membrane. For the purposes of our model, we consider the cortex to represent the combined membrane and cortex.
Actin protrusions from the cortex are allowed to bind to ECM fibers. Upon binding, the cortex stiffens, which allows the cell to ``pull'' on the ECM by generating traction forces on the tip of the protrusion \cite{legant2010measurement,zhu2016comparison}. Here we develop a model of this mechanism in 2D to gain insight into how ECM stiffness affects the ability of the cell to migrate before developing a computationally expensive 3D model.

We consider a 2D cross-section of a cell migrating through an ECM consisting of fibers immersed in fluid. The cell and nucleus are modeled as thin 1D elastic boundaries. The ECM consists of long thin fibers in 3D, and we model the cross section of one fiber as a regularized point force in our 2D model. For the cortex and nucleus, fiber elasticity gives the force density (in pN/$\mu$m$^2$) on a given configuration by
\begin{equation}
\label{eq:fibelforce}
    \bm{F}^{el}_{n/c}=\frac{\partial}{\partial s}\left(T_{n/c}\bm{\tau}\right), 
\end{equation}
where $n/c$ stands for the nucleus or cortex, $s$ is the reference arclength variable,  $\bm{\tau}=\bm{X}_s/\norm{\bm{X}_s}$ is the unit tangent vector, and 
\begin{equation}
\label{eq:fibeltension}
    T_{n/c}=k_{n/c}\left(\norm{\bm{X}_s}-1\right)
\end{equation}
is the fiber tension. $k_{n/c}$ (pN/$\mu$m) represents the stiffness of the nucleus/cortex. At the beginning of our simulations, we choose the cortex to be relatively soft with $k_c=1$ pN/$\mu$m and the nuclear boundary to be much stiffer, $k_n=50$ pN/$\mu$m \cite{friedl2011nuclear}. We take the diameter of the cortex to be 1 $\mu$m and the diameter of the ``nucleus'' to be 0.9 $\mu$m, with the latter taken to be large to model effective elasticity of the cytoplasm. 

We also discretize the cortex with $N_c=80$ points and the nuclear boundary with $N_m=40$ points. Using this discretization, one can numerically approximate derivatives in Eqs.\ \eqref{eq:fibelforce} and \eqref{eq:fibeltension} via centered differences to obtain a force density at each point in pN/$\mu$m$^2$, then multiply by the reference point spacing to obtain a force, $\hat{\bm{F}}^{el}_{n/c}$ in pN/$\mu$m at each point on the nucleus/cortex. 

The ECM can be represented in a cross sectional sense as an array of points in space with some characteristic length spacing. In this section, we keep the spacing constant and fix it to be on average the same as the diameter of the cell, so that the cell is not sterically hindered from passing through the ECM. Future work will focus on the effect of ECM density in a more rigorous context; our goal here is instead to show the effect of matrix stiffness at constant fiber density. 

\begin{figure}
\label{fig:motileexample}
\centering     
\subfigure[ECM architecture]{\label{fig:ecmtri}\includegraphics[width=60mm]{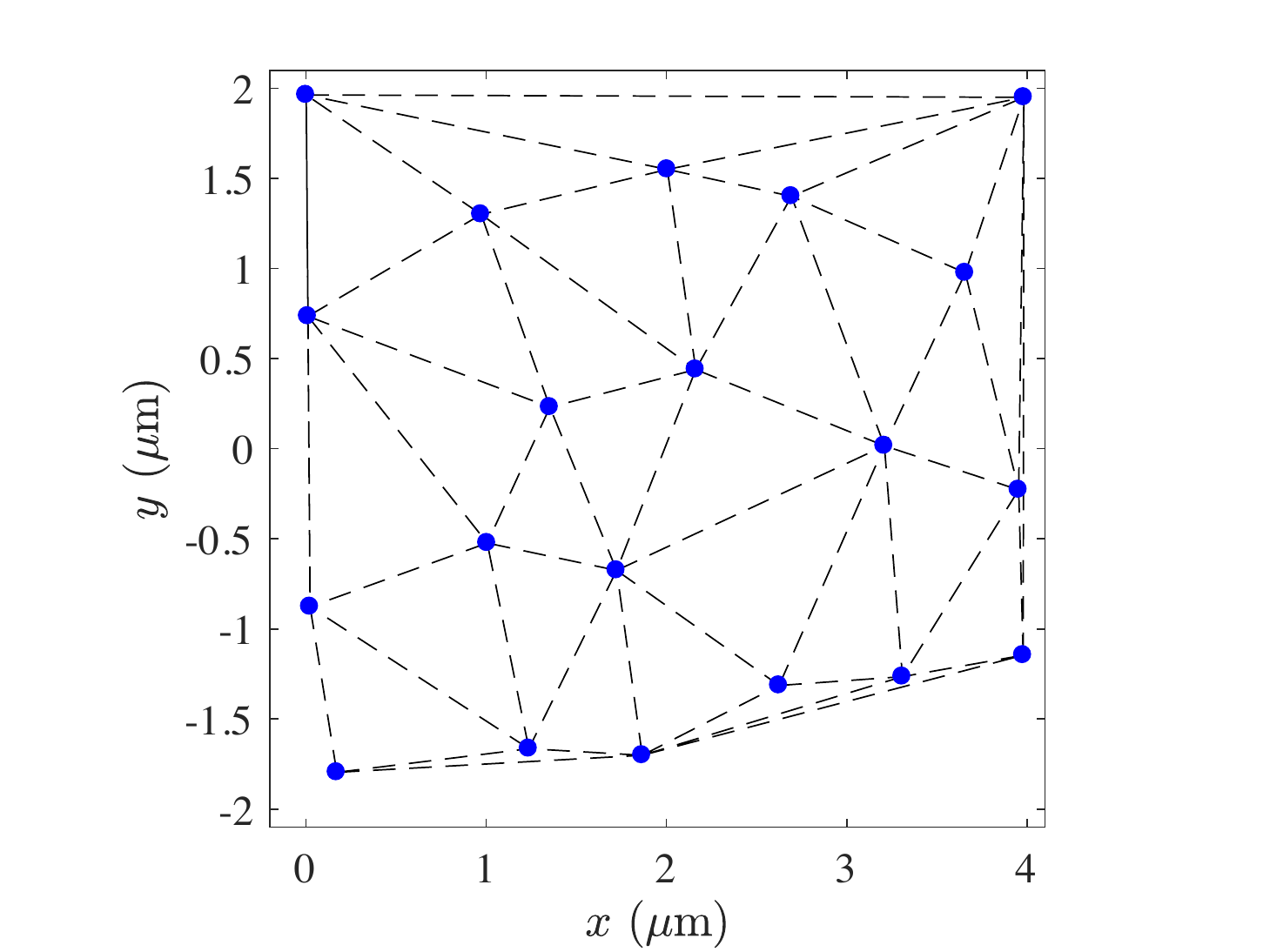}}
\subfigure[Protrusion formation]{\label{fig:mech1t0}\includegraphics[width=60mm]{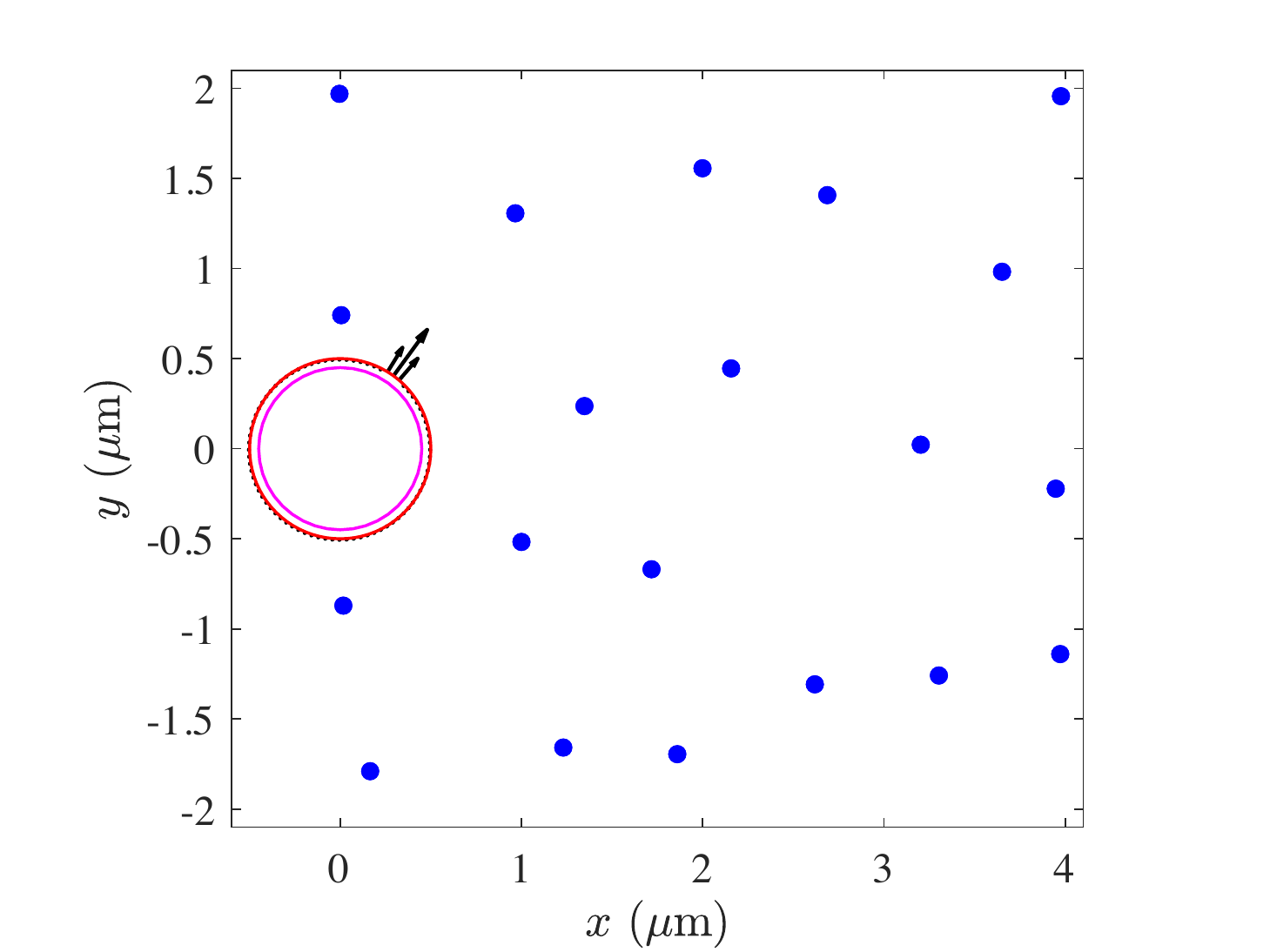}}
\subfigure[Attachment to ECM node]{\label{fig:mech1t1}\includegraphics[width=60mm]{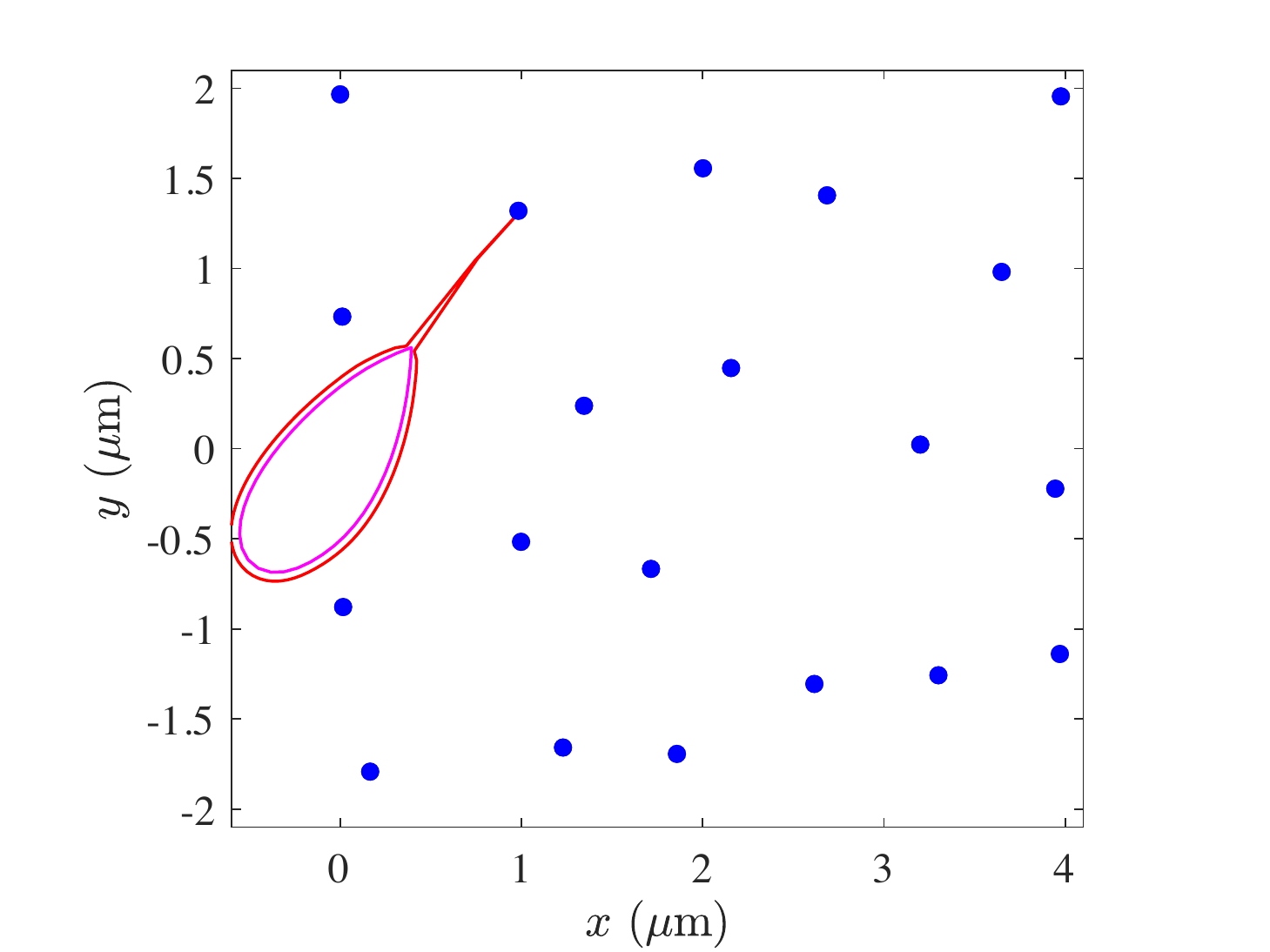}}
\subfigure[Cortical contraction]{\label{fig:mech1t2}\includegraphics[width=60mm]{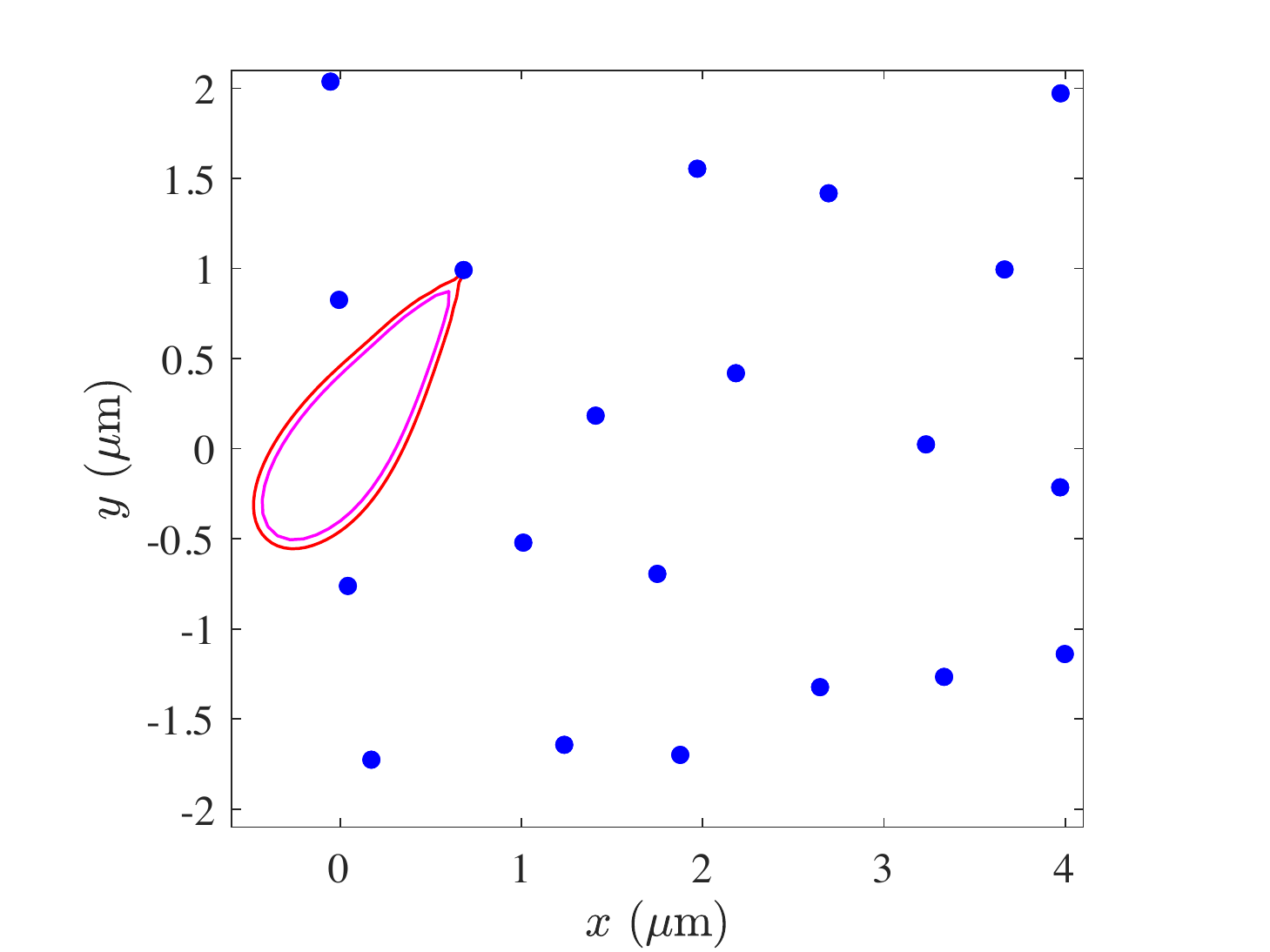}}
\subfigure[ECM pull]{\label{fig:mech1t3}\includegraphics[width=60mm]{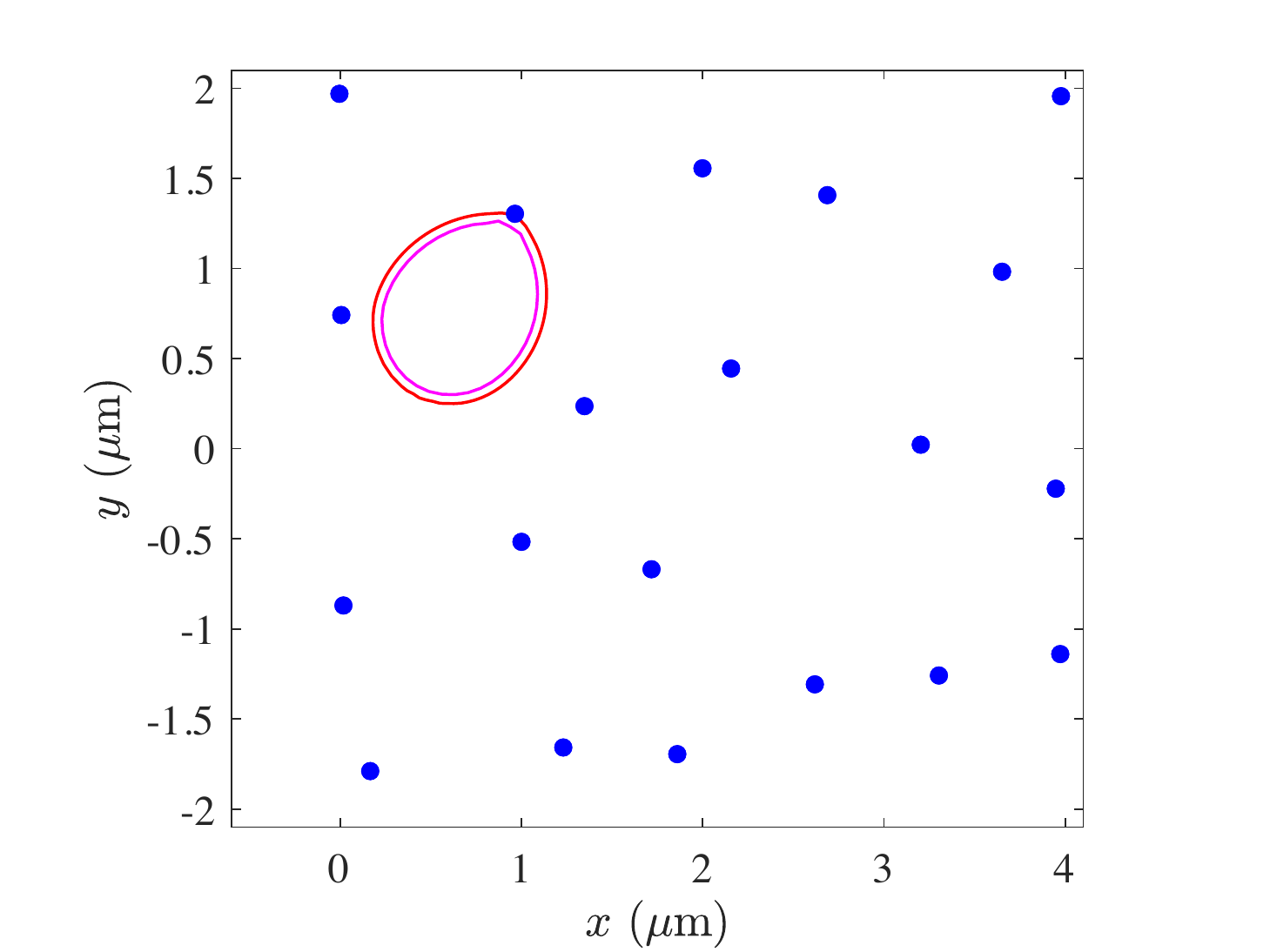}}
\subfigure[Final state]{\label{fig:mech1t4}\includegraphics[width=60mm]{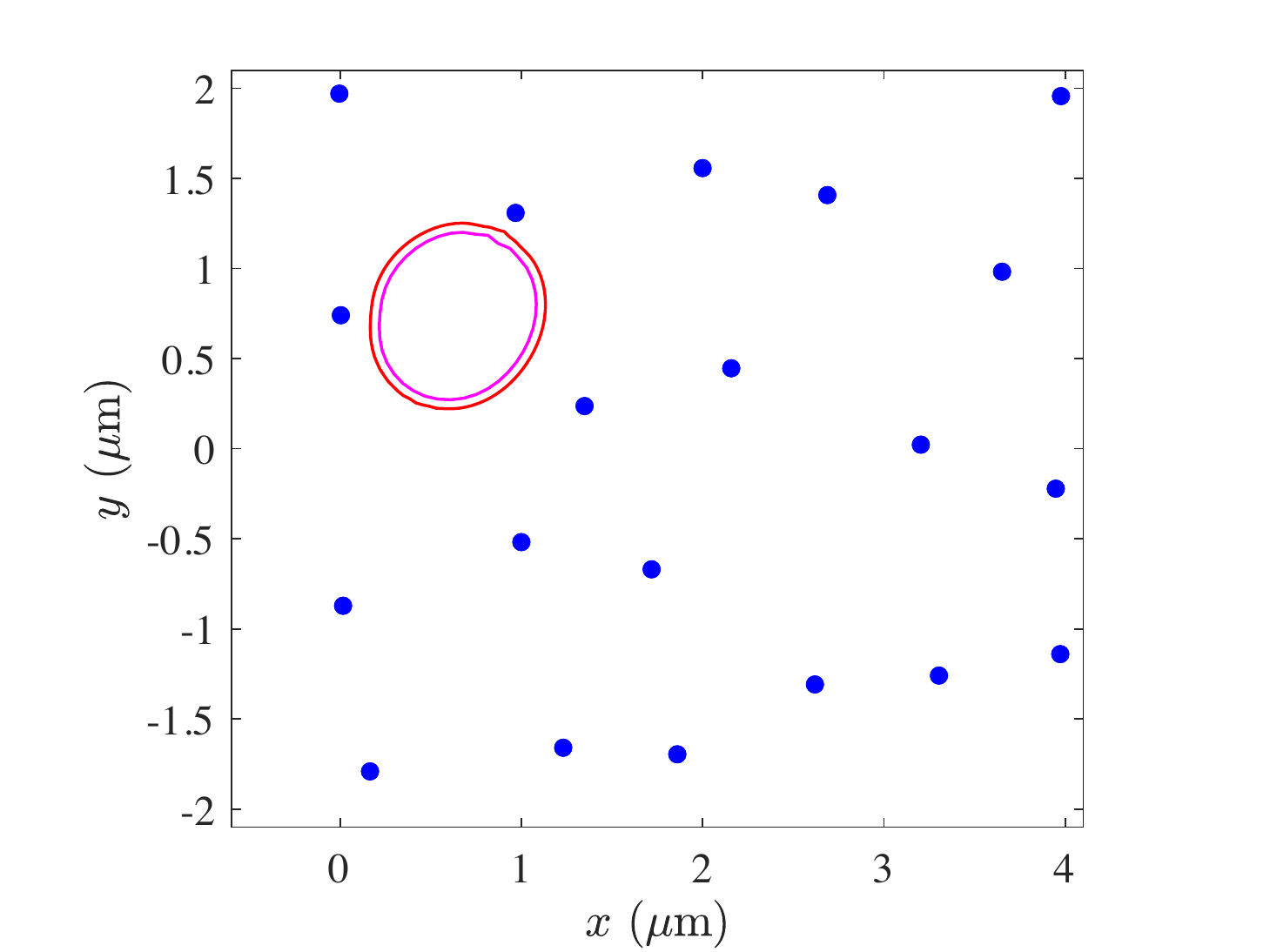}}
\caption{One cycle of a cell migrating via a finger protrusion mechanism through an ECM matrix of elastic nodes with $k_{teth}=50$ pN/$\mu$m$^2$. (a) The structure of the ECM, which has 20 nodes (blue points) that are linked together by springs (dashed black lines). (b-f) The dynamic process of cell migration. (b) A protrusion forms on the cell surface. (c) The protrusion binds to a node. (d) The cortical stiffness increases, pulling the node inward. (e) The dynamic balance between elasticity of the cell and ECM elasticity pulls the cell towards the ECM node's resting position. (f) The cell releases the node and is ready to form another protrusion.}
\end{figure}

We therefore generate 20 ECM nodes that are approximately spaced by the cell diameter on a $4 \times 4$ box, shown as blue points Fig.\ \ref{fig:ecmtri}. We triangulate this set of points, with each edge representing a spring that connects two ECM nodes (dashed black lines in Fig. \ref{fig:ecmtri}). Let $k_{teth}$ (pN/$\mu$m$^2$) denote the stiffness of these springs. Then the time-dependent force (in pN/$\mu$m) on each ECM node is given by
\begin{equation}
\label{eq:fecm}
\hat{\bm{F}}^j_{ECM}(t)= -k_{teth}\left(\bm{X}^j(t)-\bm{Z}^j+ \sum_{i \in \mathcal{N}(j)} (\bm{X}^j(t)-\bm{X}^i (t))\right). 
\end{equation}
Here $i \in \mathcal{N}(j)$ denotes a point $i$ which is a neighbor of point $j$, in the sense that the nodes are connected by an edge in the Delanuay triangulation (black dotted lines in Fig.\ \ref{fig:ecmtri}). We note also the presence of an anchoring (tether) node, $\bm{Z}^j$, whose purpose is to make sure the network stays in place dynamically. Without linking the nodes to reference nodes $\bm{Z}^j$, the force function in Eq.\ \eqref{eq:fecm} would be translation-invariant, and the entire network of nodes would be free to slide away from the cell without penalty. We can compute $\bm{Z}^j$ for each node by setting the force at $t=0$ in Eq.\ \eqref{eq:fecm} to zero and solving for $\bm{Z}^j$. This is desired physically for the cell to migrate relative to the ECM. We have therefore determined the forces on the nucleus, cortex, and ECM that need to be computed at each time point and passed to Eq.\ \eqref{eq:totalvel}. 

We note that the use of points for the ECM fibers necessitates the use of a regularized method, as the velocity due to a point force is technically infinite at that point. We set $\epsilon=0.075$ $\mu$m in the regularized equations, so that each point has an effective radius around it that is much smaller than the radius of the cell. We note that this value for $\epsilon$ is also the approximate spacing between the discrete nuclear and cortical points, which is one criterion for choosing $\epsilon$ \cite{cortez2001method, cortez2005method}. 

The overall simulation algorithm is as follows. Draw from a uniform distribution a point $j$ on the front edge of the cortex (representing cell polarization) and suppose that actin polymerizes at that point. We apply a force density of strength $f_0=500$ pN/$\mu$m$^2$ in the normal direction at point $j$ and a force density of strength $f_0/2$ in the normal direction at points $j-1$ and $j+1$. Importantly, the cell physically cannot generate any net force on the fluid, so we spread the equal and opposite force over the other $N_c-3$ cortex points. The force distribution on the cortex is shown in Fig.\ \ref{fig:mech1t0}. We note that the only effect of $f_0$ is to set the timescale of migration, and  we are concerned with the \textit{relative} timescale across different ECM stiffnesses (so that the choice of $f_0$ is arbitrary). For this reason we also set $\mu=1$ Pa-s for simplicity. 

We then allow the cell protrusion to grow by evolving the system in time until the protrusion tip comes into contact with a node. By contact, we mean that the discrete points are a distance $2\epsilon$ or less from each other, so that their ``blob'' functions are in contact. Once the discrete points come within a distance $2\epsilon$ of each other, the protrusion tip binds to the node (shown in Fig.\ \ref{fig:mech1t1}), and the cortex becomes stiffer by a factor of 100 to model the increased traction at the protrusion tips seen in \cite{legant2010measurement}. The increased stiffness causes the cortex to rapidly become rounder. Since the cortex is attached to the node, it then pulls the node inward as shown in Fig.\ \ref{fig:mech1t2}. As the node is pulled in, it generates a force in the opposite direction due to elasticity of the ECM (the node's resting configuration is its initial configuration in the ECM lattice). These forces balance dynamically, so that as the cortex becomes more round, the cortical force due to elasticity decreases, which in turn allows the force on the ECM to decrease, thereby pulling the entire cell and node back towards the initial position of the node. In the final state, Fig.\ \ref{fig:mech1t3}, the cell is round and the node returns to a point close to its initial position. At this time, the node detaches from the cell by moving a distance $2\epsilon$ away in the normal direction, as shown in Fig.\ \ref{fig:mech1t4}, and the process can then repeat. We define this entire process as one cycle.   

\begin{figure}
\centering     
\subfigure[$R=10^3$]{\label{fig:ourfixfinal}\includegraphics[width=60mm]{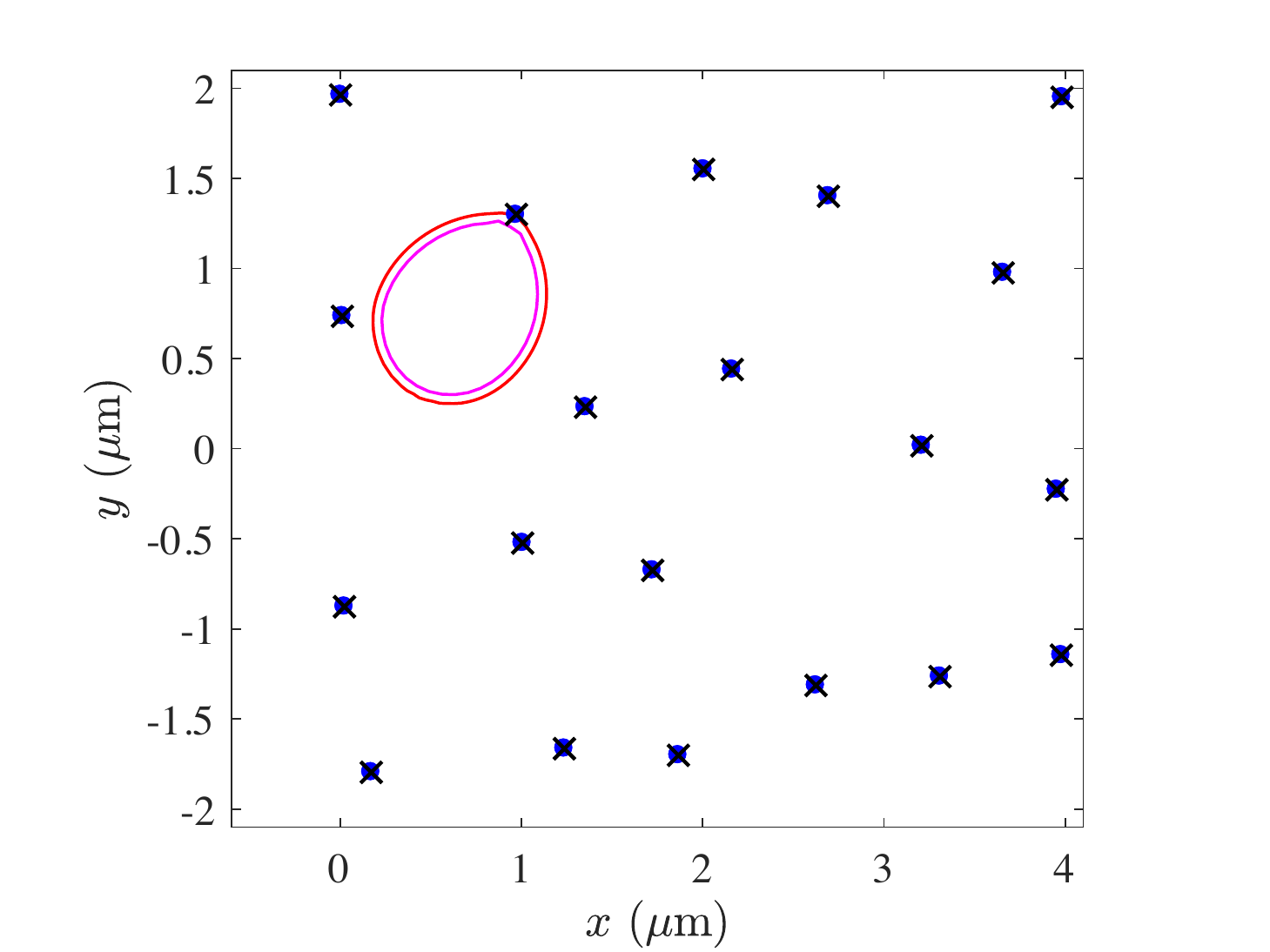}}
\subfigure[Subtracting force mean]{\label{fig:meansubfinal}\includegraphics[width=60mm]{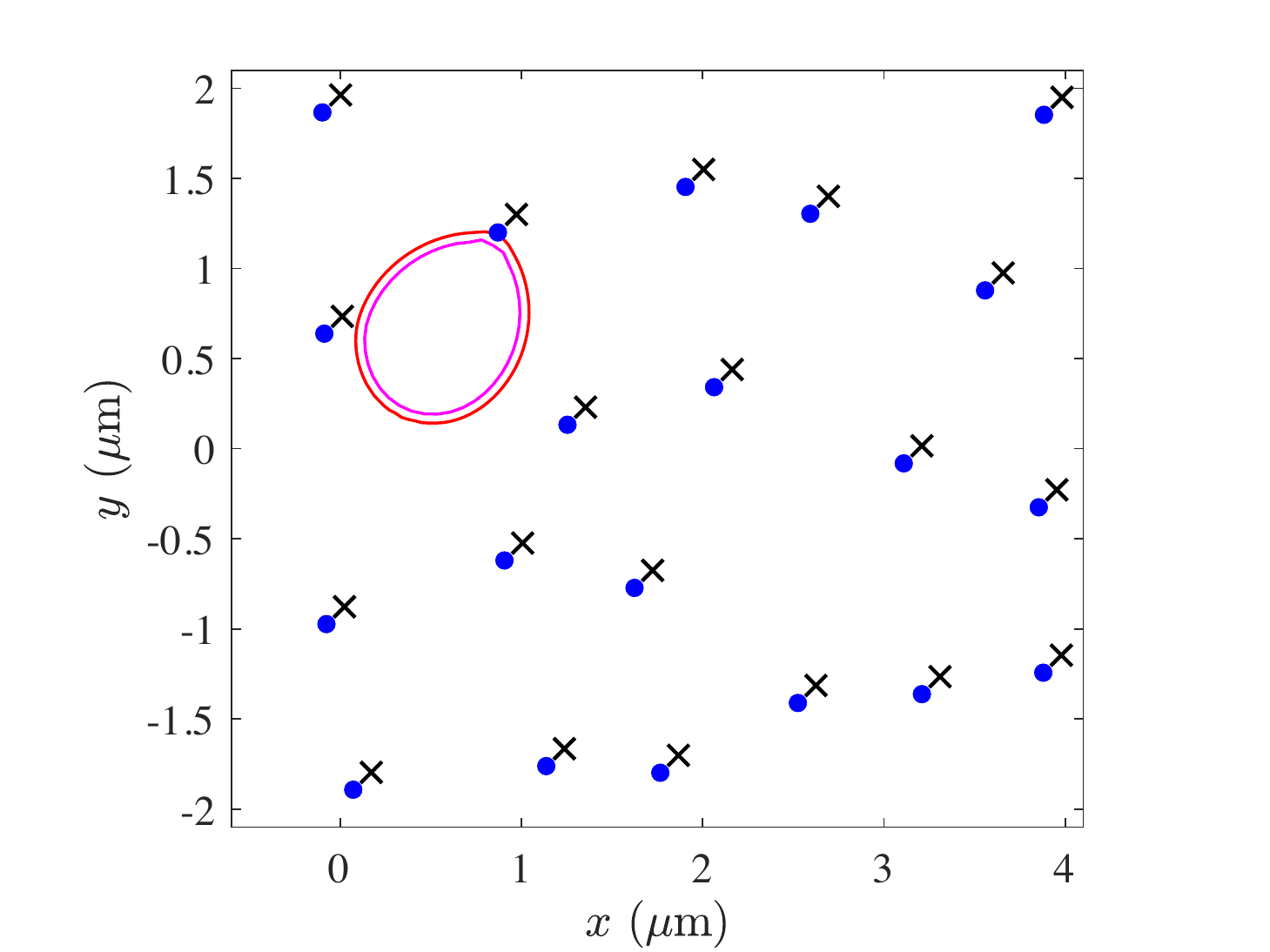}}
\caption{Final states for one cycle of the cell motility problem with $k_{teth}=50$ pN/$\mu$m$^2$ and (a) $R=10^3$ (b) sum of forces being zero via mean subtraction. Subtracting the mean of the force has created an artificial equilibrium state in (b), where each of the nodes has been displaced from its original position (black x's) by some constant amount. This does not occur in (a).}
\label{fig:meanRcomparemot}
\end{figure}

In this application, the anchor ECM nodes, $\bm{Z}^j$ create a net force in the domain. As we observed in Section \ref{sec:motex}, handling this imbalance by subtracting the mean force at each node can create non-physical translated equilibrium configurations. Fig.\ \ref{fig:meanRcomparemot} shows the final configuration when the system velocity has dropped below $\epsilon$ for a migrating cell in an ECM with $k_{teth}=50$ pN/$\mu$m$^2$ for (a) a system simulated using Eq.\ \eqref{eq:totalvel} and (b) a system simulated with zero net forcing via subtracting the mean force at each node. A shift in the entire domain in the negative horizontal and vertical directions is shown in Fig.\ \ref{fig:meansubfinal}.  Physically, we expect the nodes to return to their initial configuration in Fig.\ \ref{fig:ecmtri} (marked with black x's in Fig.\ \ref{fig:meanRcomparemot}). 

The translation of the ECM structure in the case of subtracting the mean forces occurs because the pulling inward of the ECM node (shown in Fig.\ \ref{fig:mech1t2}) creates a net force in the positive horizontal and vertical directions. Subtracting the mean force from each node does result in a net applied force of zero, but also results in an equilibrium configuration where the nodes have been shifted in the direction opposite the force imbalance. This situation is analogous to that of Section \ref{sec:motex}, where there was an artificial equilibrium state in the positive horizontal direction (and no relaxation to the equilibrium) resulting from a force imbalance in the negative horizontal direction. Such a shift may not be important for some applications if the \textit{relative} position of the objects is desired. For our application, we are interested in the \textit{absolute distance traveled by the cell}. For this reason, along with the ECM returning to its initial resting configuration for subsequent motility cycles, we conclude that it is better to use Eq.\ \eqref{eq:totalvel} to update the velocity rather than subtracting the mean forces.

\begin{figure}
\centering     
\includegraphics[width=60mm]{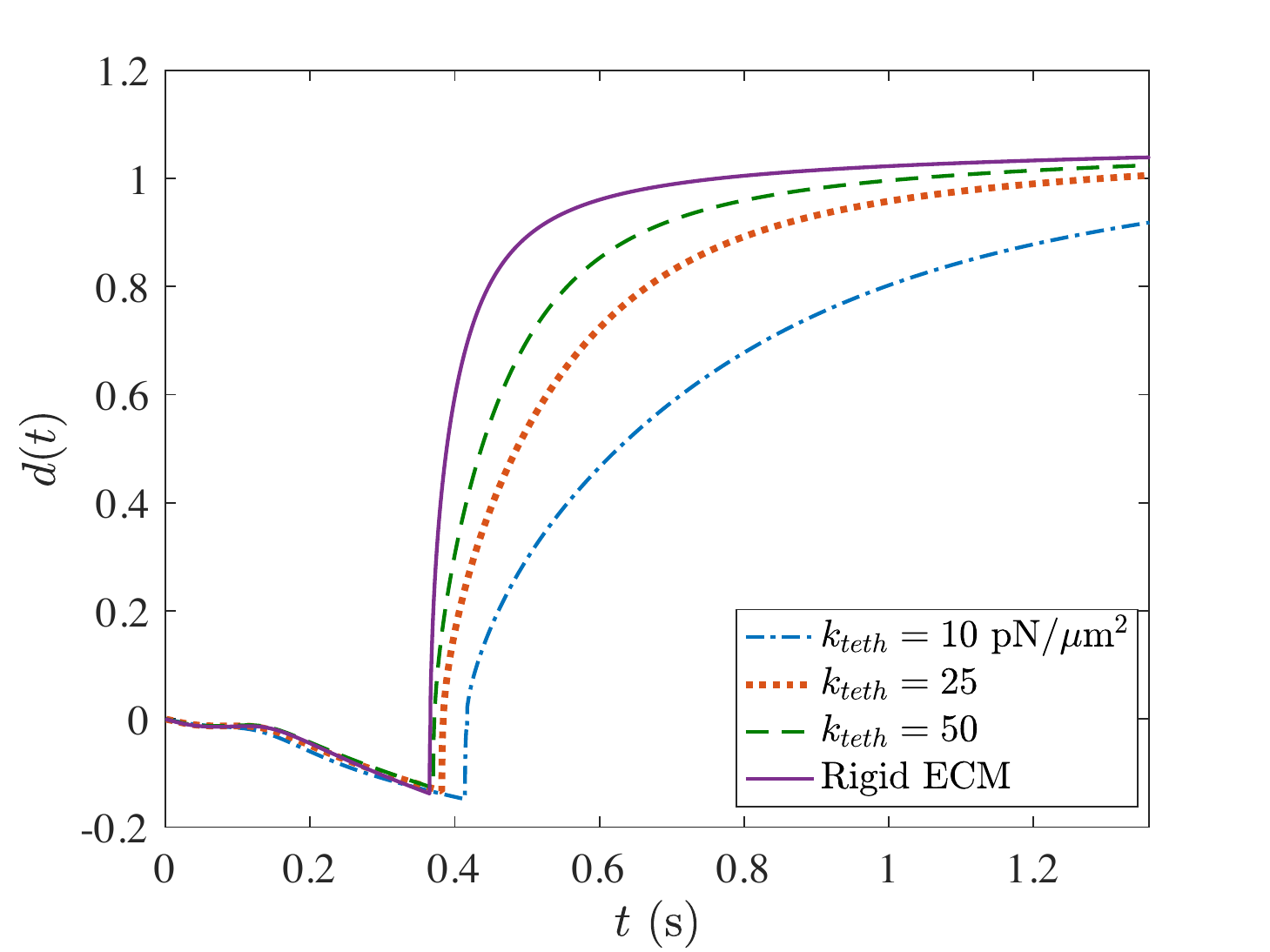}
\caption{Euclidean distance traveled by the nucleus' center of mass vs. time for different values of the ECM stifness. Stiffer ECMs display faster velocities.}
\label{fig:motileks}
\end{figure}

We now use the model to simulate the distance traveled by the cell for different values of ECM stiffness $k_{teth}$. We calculate the total Euclidean distance traveled by the nucleus' center of mass as a function of time for the same 20 node ECM, but with varying stiffness $k_{teth}=10, 25, 50$ pN/$\mu$m$^2$. For comparison, we also simulate a rigid ECM by enforcing a $\bm{u}=0$ boundary condition at each of the ECM nodes rather than the mean velocity condition on the large circle. We simulate up to a finite time, which corresponds to the time the cell has finished one cycle of migration (Fig.\ \ref{fig:mech1t3}) in the rigid ECM case ($t=1.36\approx1360\Delta t$). The timestep is adaptive; generally it is taken to be $\Delta t=0.001$, but it shrinks to $\Delta t=2 \times 10^{-4}$ for a small time ($0.05$ s) beginning when the cortex binds to a node and stiffens to $k_c=100$ pN/$\mu$m. In Fig.\ \ref{fig:motileks}, we plot the Euclidean displacement in the direction of the ECM node over time for different values of stiffness $k_{teth}$. In all cases, the cell initially moves backwards slightly (as seen in Fig.\ \ref{fig:mech1t1}) prior to contacting an ECM node. Once the ECM node is contacted and the cortex contracts, the distance traveled increases with time, with larger velocities for stiffer matrices. However, the data from all simulations appear to be approaching the same steady state value of displacement. For stiffer matrices, the node resists deformation by the cell (shown in in Fig.\ \ref{fig:mech1t2}), and the cell moves toward the node. In the rigid case, the node does not move, and the cell quickly contracts to form a circular configuration around the node. 
The conclusion of this preliminary study is therefore that \textit{stiffer matrices} allow for faster cell velocities for finger-like protrusion mechanisms. We plan to study this problem in more detail by varying the matrix density and nuclear and cortical stiffness in future work.

\subsection{Cellular blebbing}
\label{sec:bleb}
Cellular blebs are spherical membrane protrusions that have been observed during cell migration \cite{CharrasPaluch}. A  bleb forms when the cell cortex, normally attached to the cell membrane by linker proteins, detaches from the membrane. Cells that bleb are pressurized due to actomyosin contractility within the cortex. Once a bleb is initiated, a pressure driven flow drives the intracellular fluid (cytoplasm) that locally expands the membrane. 

\begin{figure}
\centering     
{\includegraphics[width=60mm]{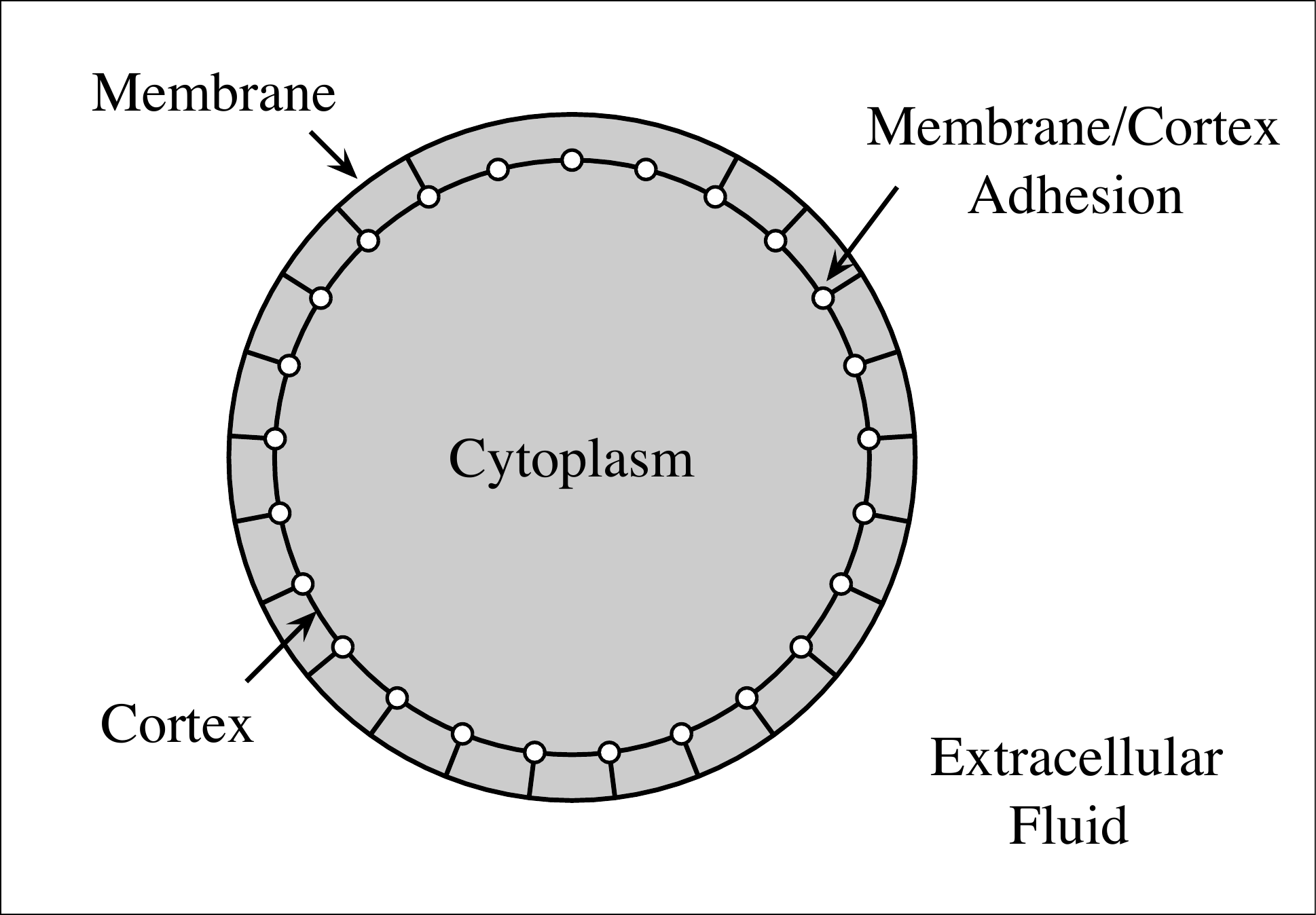}}
\caption{Components of the bleb model. The cytoplasm is modeled as a viscous fluid. A bleb is initiated by removing membrane-cortex adhesive links in a small region at the top of the cell.}
\label{fig:blebscheme}
\end{figure}

Bleb initiation has been modeled using different approaches, including solid mechanics \cite{woolley2014}, the immersed boundary (IB) method \cite{strychalski2013computational,strychalski2016intracellular}, and boundary integral methods \cite{fang2017combined,lim2012size,lim2013computational}. Results from several of these models 
have shown a bleb relieves only a small amount of intracellular pressure when the cytoplasm is modeled as a viscous fluid \cite{strychalski2013computational,tinevez2009role}. 
Results from other models simulated with boundary integral methods show large pressure relief after bleb expansion \cite{fang2017combined,lim2013computational}.
Our goal is to identify the source of this contradiction because maintaining high intracellular pressure is essential for cells to migrate using blebs \cite{paluchraz2013}.

Here we present a model of bleb expansion based on \cite{lim2013computational}. 
We treat both the membrane and cortex as one dimensional closed curves. The membrane and cortex parameterizations 
are represented by $\bm{X}^m(s)$ and $\bm{X}^c(s)$, respectively, where $s$ is the arclength parameter. The most critical part of the model is the adhesion that connects the membrane and cortex. We model adhesion by an elastic spring connecting the membrane to the cortex with stiffness $k_{adh}$. The force density on the membrane due to adhesion is given by, 
\begin{equation}
\label{eq:Fmemcor}
\bm{F}^{mem/cor}_{adh}(s) = -k_{adh}\left(\bm{X}^m(s) - \bm{X}^c(s) \right), 
\end{equation}
with the force density on the cortex, $\bm{F}^{cor/mem}_{adh}(s)$,  equal and opposite.
Elastic forces on the membrane and cortex are due to surface tension and stretching and are computed by Eq.\ \eqref{eq:fibelforce} with
\begin{equation}
T=\gamma_m+k_m(\norm{\bm{X}_s}-1)
\end{equation}
with constants $\gamma_m$ and $k_{m}$ representing membrane surface tension and stiffness, respectively. The corresponding elastic parameters for the cortex are denoted by $\gamma_c$ and $k_c$.
The membrane satisfies a no slip boundary condition, and its velocity is computed by Eq.\ \eqref{eq:totalvel}.
The velocity of the cortex is computed via a force balance, similar to \cite{lim2013computational}, 
\begin{equation}
\label{eq:cormov}
\frac{d\bm{X}^c}{dt} = \frac{1}{\nu_c}\left(\bm{F}^{cor}_{el}+\bm{F}^{cor/mem}_{adh} \right),
\end{equation}
where $\nu_c$ is the cortical viscosity. A bleb is initiated by removing the adhesive links in a small region of length approximately 5 $\mu$m at the top of the cell (see Fig.\ \ref{fig:blebscheme}).  


\begin{table}
\centering
\begin{tabular}{|lllc|} 
\hline
Symbol & Quantity & Value & Source \\
\hline
$r_{\rm mem}$ & Cell radius & 10 $\mu$m & \cite{strychalski2016intracellular,tinevez2009role}\\ 
$r_{\rm cortex}$ & Cortex radius & 9.9 or 9.85 $\mu$m & \cite{fang2017combined} \\
$\gamma_{m}$ & Membrane surface tension& 40 pN/$\mu$m & \cite{strychalski2016intracellular}\\
$k_m$ & Membrane stiffness & 80 pN/$\mu$m &  \\
$\gamma_c $ & Cortex surface tension & 250 pN/$\mu$m & \cite{strychalski2013computational}  \\
$k_c$ & Cortical  stiffness & 100 pN/$\mu$m & \cite{strychalski2013computational}  \\
$k_{\rm adh}^{\rm mem/cortex}$ & Membrane/cortex adhesion & $247 $ pN/$\mu$m$^3$ &   \\
& stiffness coefficient &  & \\
$\mu$ & Cytosolic viscosity & $5$ Pa\--s &  \cite{strychalski2013computational} \\
$\nu_c$ & Cortical viscosity & 10 pN-s/$\mu$m$^3$  & \cite{strychalski2013computational} \\ 
\hline
\end{tabular}
\caption{Parameters for the blebbing model.}
\label{ModelParametersNoCyto}
\end{table}


\begin{figure}
\centering     
\subfigure[$t = 0.1$ s] {\label{fig:t01size}\includegraphics[width=60mm]{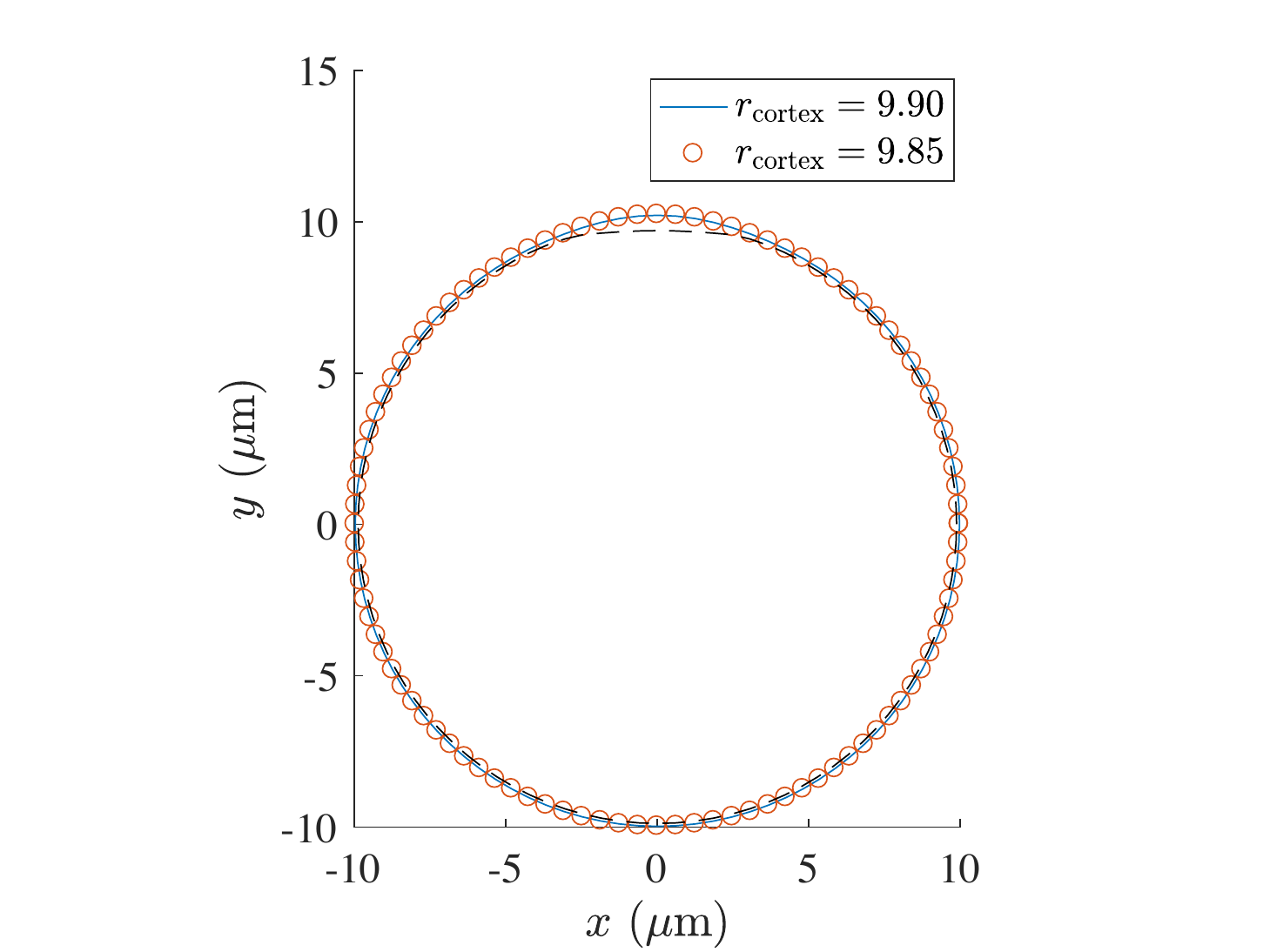}}
\subfigure[$t = 1.0$ s]{\label{fig:t1size}\includegraphics[width=60mm]{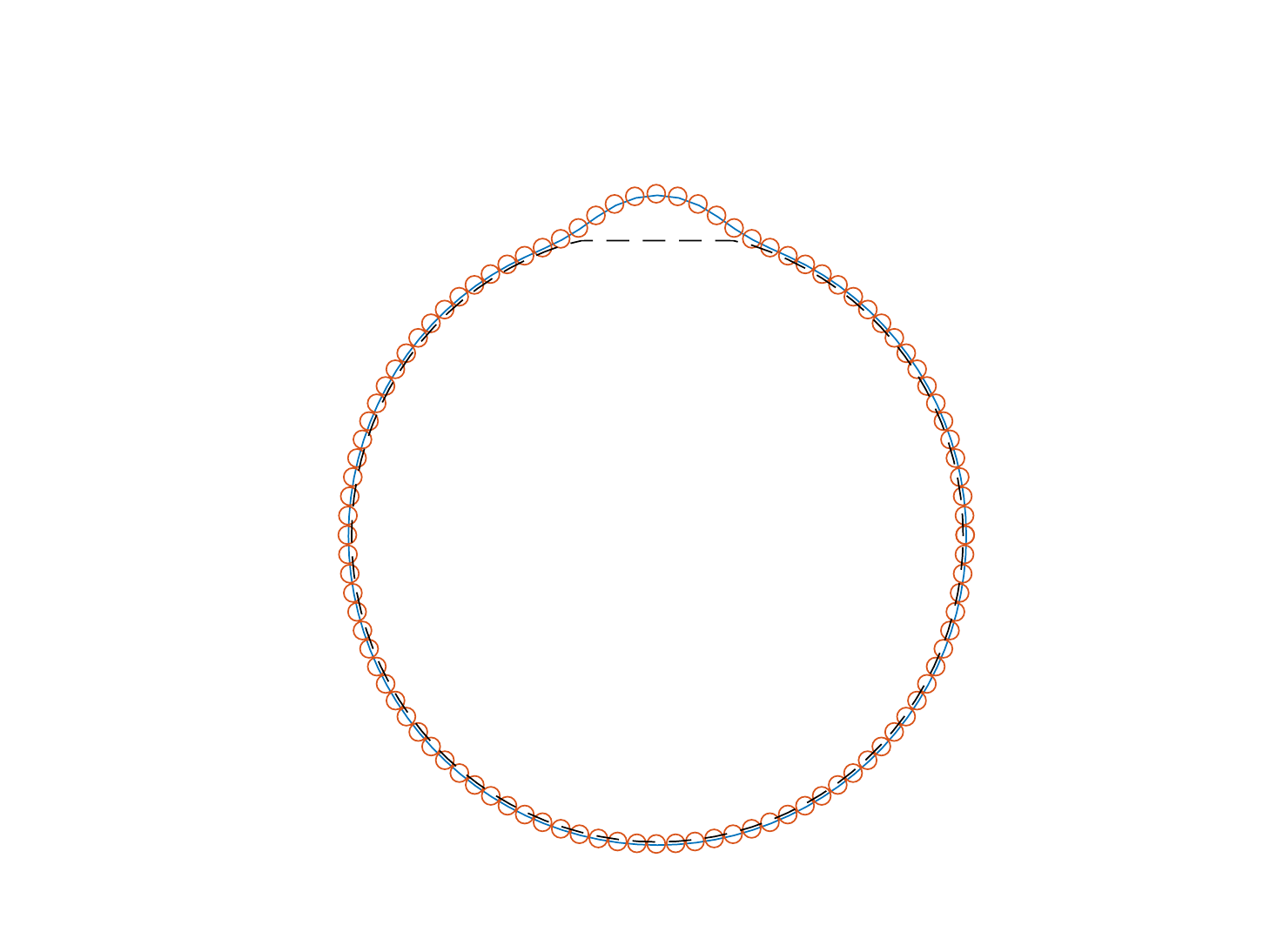}}
\subfigure[$t = 5.0$ s]{\label{fig:t5size}\includegraphics[width=60mm]{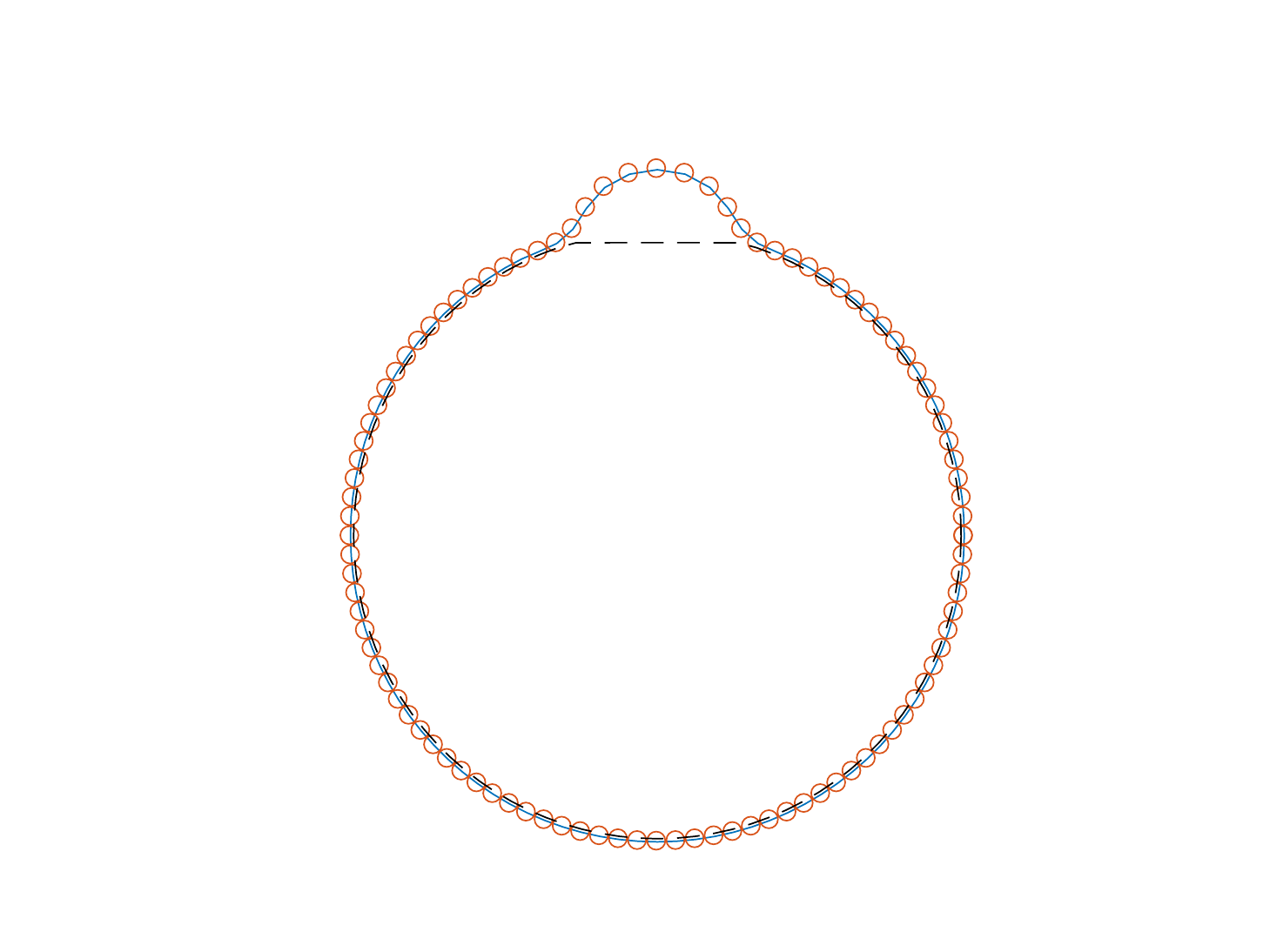}}
\subfigure[$t = 15.0$ s]{\label{fig:t15size}
\includegraphics[width=60mm]{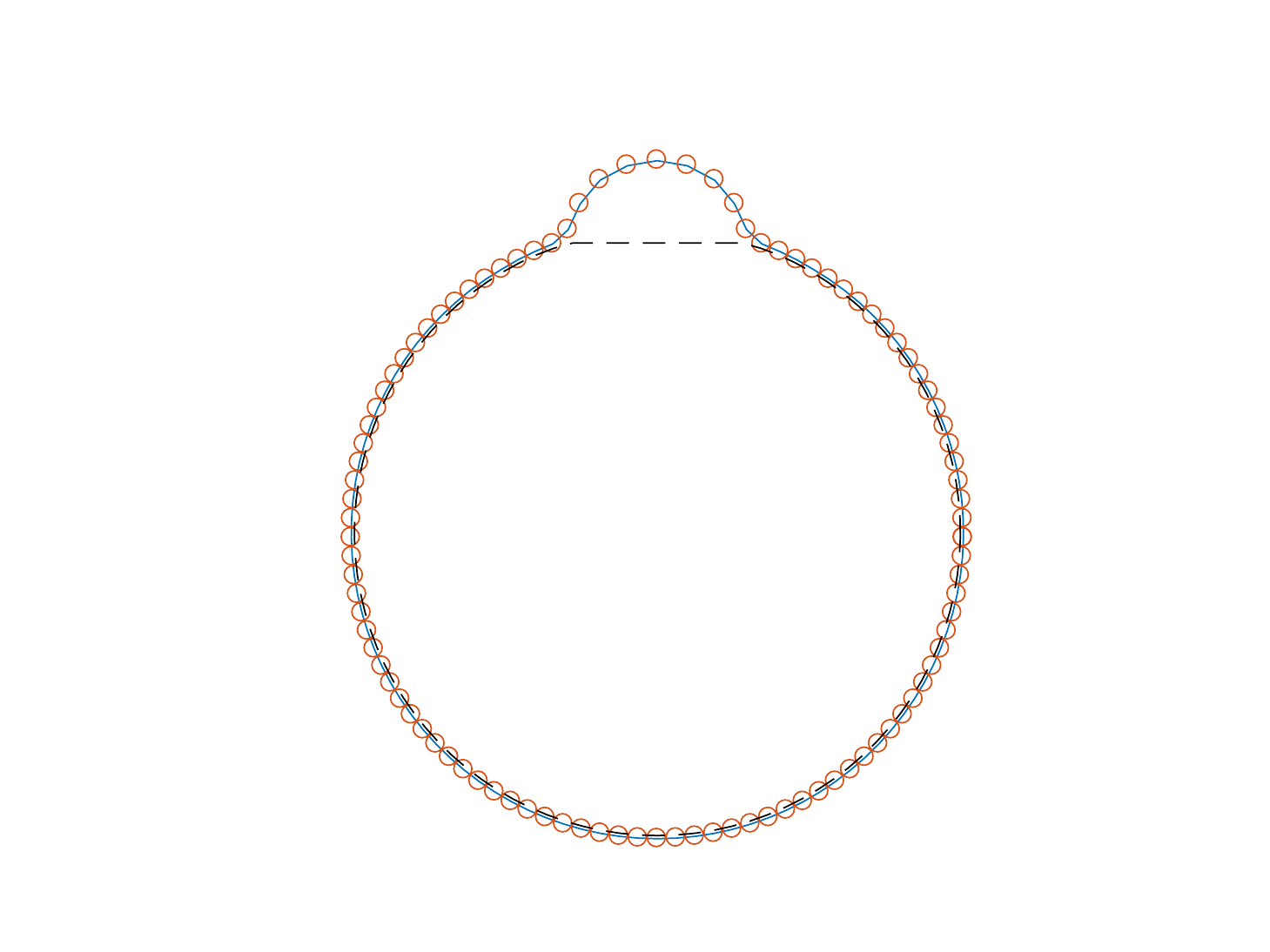}}
\caption{Membrane position for a blebbing cell with initial cortex radius 9.90 $\mu$m (blue line) or 9.85 $\mu$m (red circles). The position of the cortex is shown as a dashed black line and is in approximately the same position in both simulations. The positions are shown at several time values after bleb initiation.}
\label{fig:blebsizes}
\end{figure}

\begin{figure}
\centering     
\subfigure[$t = 0$ s] {\label{fig:t0pres}\includegraphics[width=60mm]{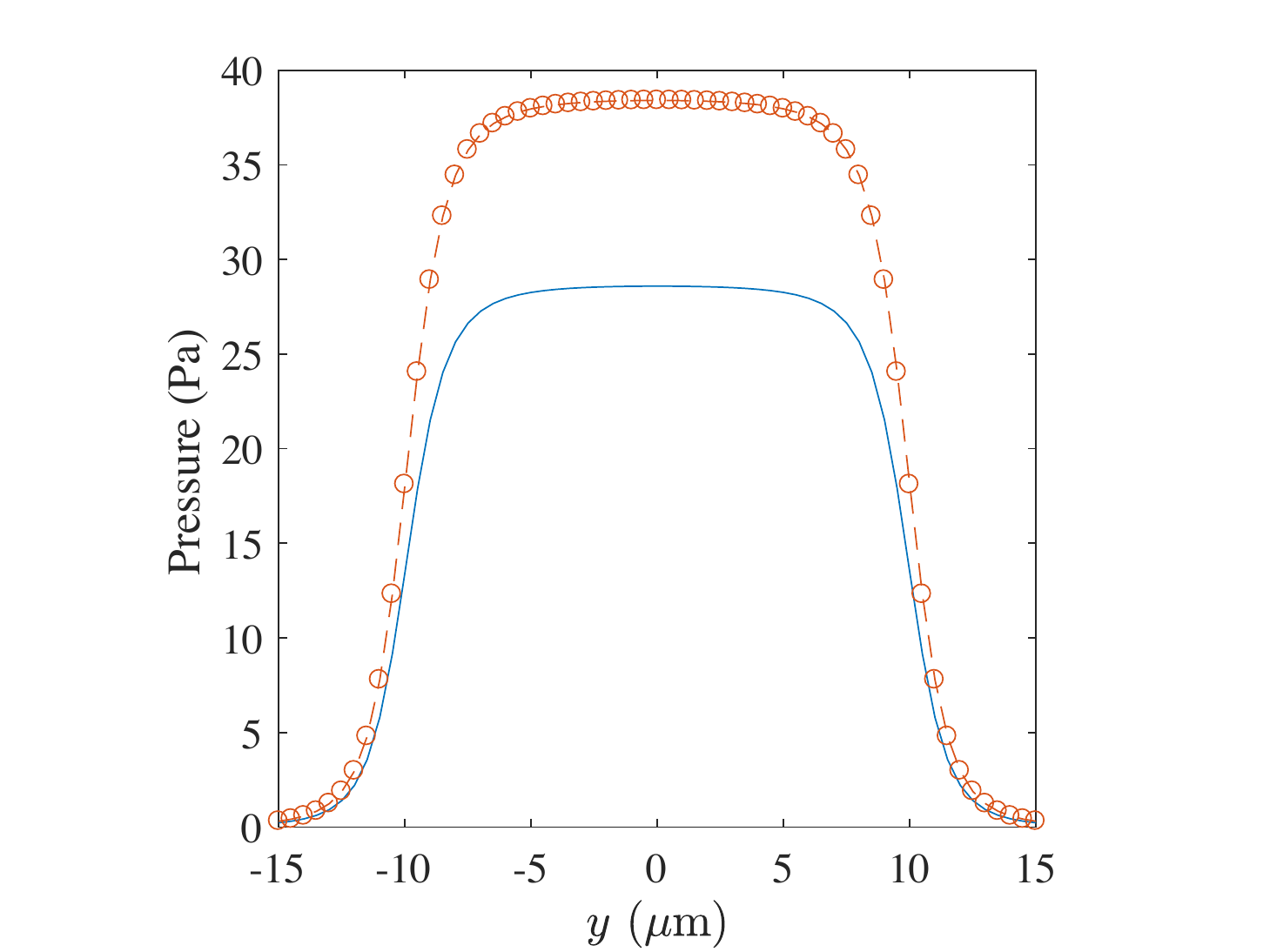}}
\subfigure[$t = 0.1$ s]{\label{fig:t01pres}\includegraphics[width=60mm]{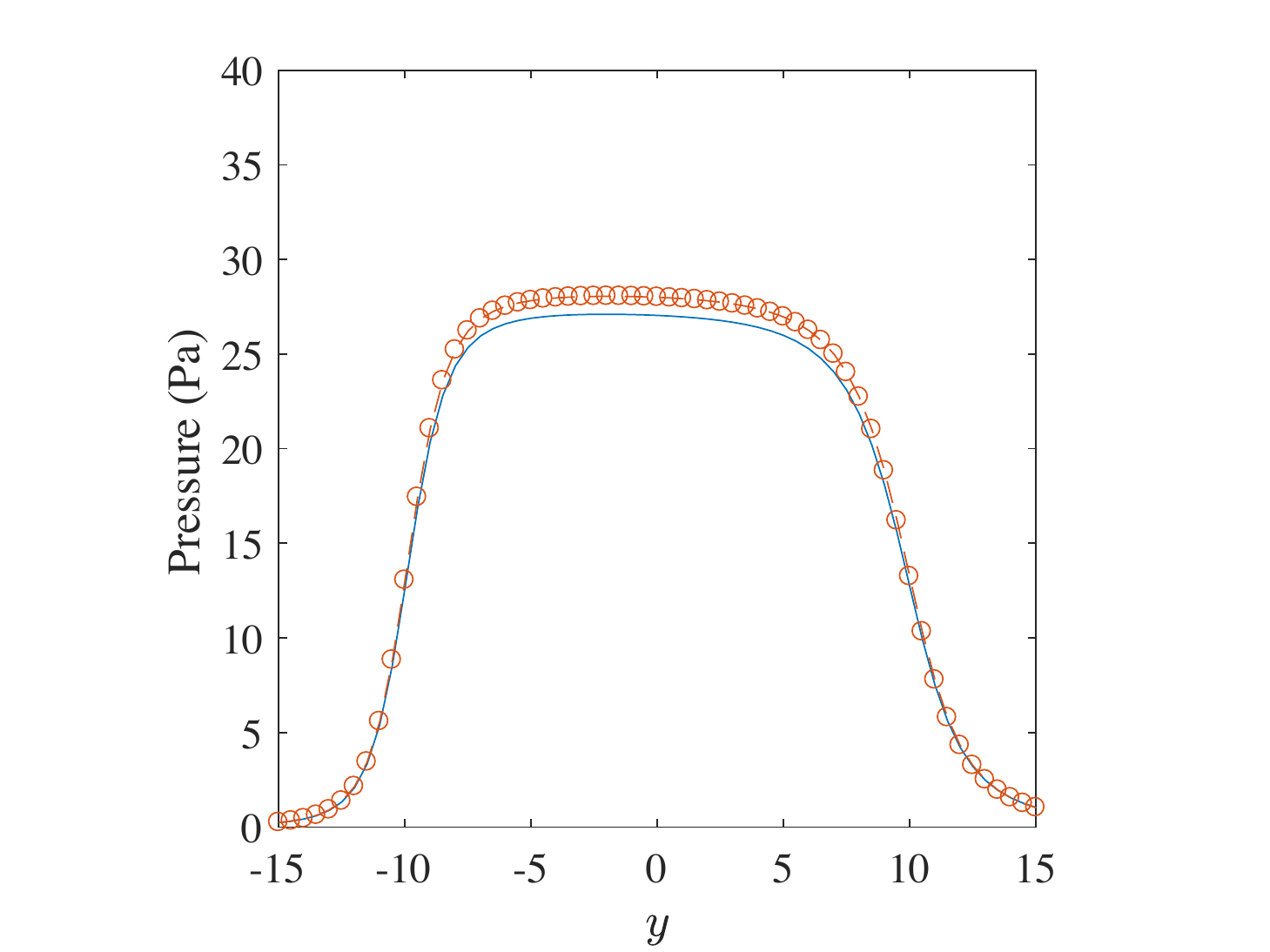}}
\subfigure[$t = 1.0$ s]{\label{fig:t1pres}\includegraphics[width=60mm]{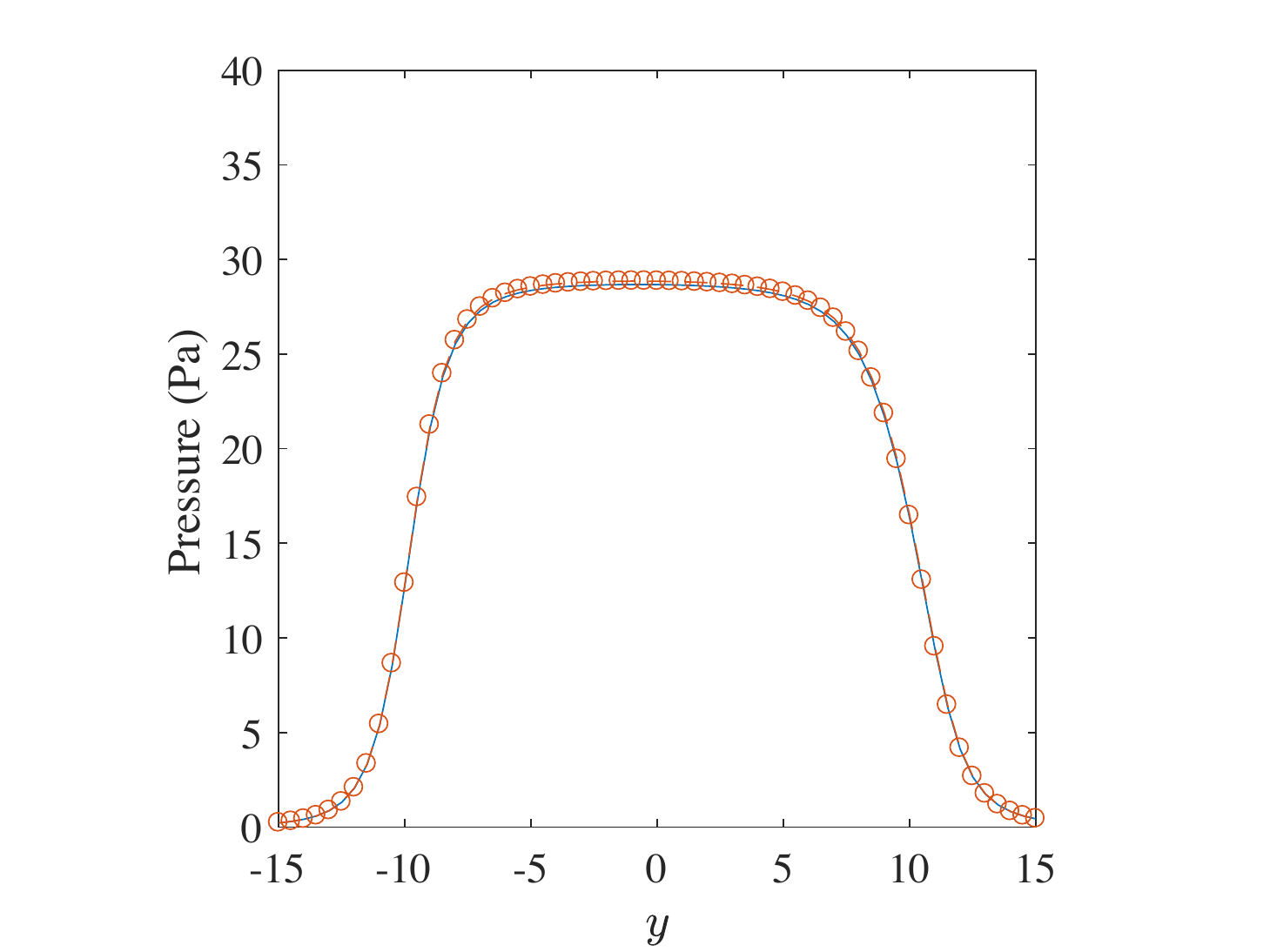}}
\subfigure[$t = 10.0$ s]{\label{fig:t10pres}\includegraphics[width=60mm]{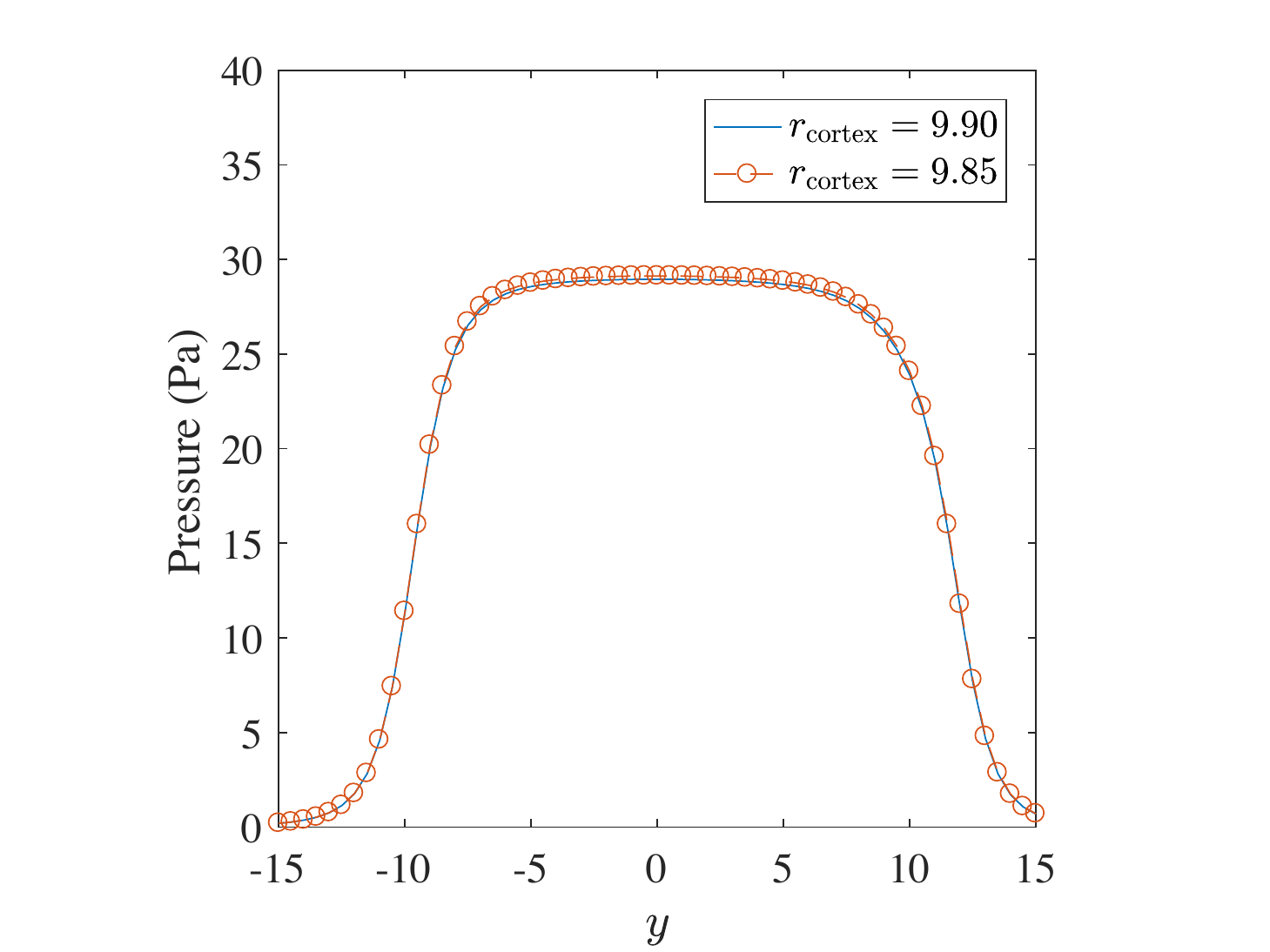}}
\caption{Pressure profile along the line $x=0$ for the blebbing cell with initial cortex radius 9.90 $\mu$m (blue line) or 9.85 $\mu$m (red circles). Profiles are shown at (a) $t=0$ s, (b) $t=0.1$ s, (c) $t=1.0$ s, and (d) $t=10.0$ s. Note the large pressure relief when the forces are initially unbalanced on the cortex ($r_{\rm cortex} = 9.85$).}
\label{fig:blebpressures}
\end{figure}

This particular blebbing model is a physical example where the net force on the membrane fiber is nonzero. When the links between the top of the membrane and cortex are broken, there is a net vertical force on the membrane because part of the adhesive forces acting in the negative vertical direction are no longer present. Even though the cortex feels the equal and opposite forces, it moves independently of the fluid according to Eq.\ \eqref{eq:cormov}. The net hydrodynamic force is therefore equal to the net force on the membrane, and is nonzero. Without including the constant velocity from Eq.\ \eqref{eq:addedvel}, the membrane would escape from the domain because of a spurious downward velocity resulting from the force imbalance in the positive vertical direction. The addition of $\bm{u}^R$ results in a well-posed problem so that we may solve for the membrane position in the blebbing model using the method of regularized Stokeslets.

Although previous studies have used the method of regularized Stokeslets to simulate cellular blebbing and migration in 2D \cite{lim2012size,lim2013computational}, the authors did not specify how they addressed the force imbalance in their models. We found our approach to give near identical results to a model where the net zero force constraint is enforced by subtracting the mean of the calculated forces at each Stokeslet point. 

We simulate the cellular blebbing process with the parameters in Table \ref{ModelParametersNoCyto}. The two different values of the cortex radius are used to test our hypothesis that a force imbalance on the \textit{cortex} is what drives the pressure relief seen by previous authors \cite{fang2017combined,lim2012size,lim2013computational}. For $k_{adh}=247 $ pN/$\mu$m$^3$, the forces on the cortex (in the absence of a bleb) are exactly in balance when $r_{\rm cortex} = 9.9 \, \mu$m. When $r_{\rm cortex} = 9.85 \, \mu$m, there is initially a force imbalance on the cortex independent of bleb initiation. 

We first equilibrate the model for ten time steps, then initiate a bleb at $t=0$ by breaking the adhesion at the 7 (out of $N=100$) points with largest $y$ coordinate. 
Fig.\ \ref{fig:blebsizes} shows the membrane shape over time for the two different values of the cortex radius, where time units are reported after bleb initiation.
The bleb sizes and shapes are exactly the same. Despite this, the pressure dynamics of the two models are quite different. As shown in Fig.\ \ref{fig:blebpressures}, the pressure drops significantly when $r_{\rm cortex} = 9.85$ but remains constant in the case $r_{\rm  cortex} = 9.90$. This is because of the force imbalance on the cortex in the former case. When the cortex's initial position is inwards of its resting position, it expands outward dynamically. This decreases the force on the membrane (and on the fluid) due to membrane-cortex adhesion (Eq.\ \eqref{eq:Fmemcor}), leading to a global pressure decrease inside the cell. 
Importantly, we observe that at $t \geq 1$ s, the cortex has reached its resting position and the two pressure profiles are the same (and are unchanged substantially with bleb expansion). 

Models that include dynamic breaking of membrane-cortex adhesive links exhibit drastic changes in intracellular pressure \cite{fang2017combined,lim2013computational}. In such models, the dynamic breaking of adhesive links over time leads to pressure relief because the membrane force is updated suddenly without accounting for the corresponding force imbalance on the cortex. As the cortex slowly responds, it contracts inward in response to the loss of adhesive force, which promotes more link breakage and pressure changes. The cortex itself is therefore never truly in equilibrium, and the force imbalance on it drives pressure changes. 

The assumption that forces on the cortex are not equilibrated may be valid during highly dynamic processes such as during cell migration \cite{YipIntBiol2014}. However, some experiments involve isolating a specific event, such as the expansion of a single bleb in \cite{tinevez2009role}. In this work, experimental data show the cell achieves a quasi-steady state behavior after bleb expansion, and the cortex is unlikely to be dynamically relieving pressure.





\section{Conclusion}
When developing models for systems from biology, physics, and engineering that involve fluid-structure interaction, the simplification from 3D to 2D allows for model prototyping and  fast simulations. In our applications, we seek to simulate cell motility and blebbing under a broad range of parameters, so fast simulations that are easy to visualize are critical for understanding model behavior. In zero-Reynolds number flow, boundary integral methods are appealing because the velocity (and position) of immersed structures can be easily computed at the locations of interest rather than by interpolation after solving for the velocity on an Eulerian grid as in the IB method \cite{peskin2002immersed}. The condition of net zero force for 2D boundary integral methods can be a limiting factor during the development and simulation of mathematical models, especially those that include elastic tether-like forcing. 

Here we present a numerical method to treat force constraints in 2D Stokes flow in an infinite domain. When the regularized or standard Stokeslet is used with a net nonzero hydrodynamic force in the flow domain, the velocity is unbounded at infinity (Stokes' paradox). For problems where no specific boundary conditions are imposed from the physics of the model system, the standard free space Stokeslet solution without modification fails because it is a Green's function for a system of equations that is not well-posed. The treatment of this problem must therefore involve solving a new, well-posed system with the appropriate Green's function.  

One option is to require the net force to be zero within the domain by subtracting the mean force, thereby making the free space problem well-posed. This approach maintains the system time and length scales, but we show here that it can lead to falsely translated equilibrium states and unphysical dynamics. Alternatively, shifting the boundary conditions to create a well-posed system and locally valid solution is appealing (previously treated via solving a linear system \cite{Copos2018,cortez2001method,Vanderlei2011}). We show here that by using a confined geometry and enforcing the condition of a mean zero velocity on the boundary, we can easily derive a new Green's function for a well-posed system of equations that gives the physically correct behavior. Introducing the new boundary has the effect of introducing a new length scale and a corresponding change in the Reynolds number, and we use this fact to derive an upper bound on the size of the boundary. We combine this with a lower bound that comes from the variation of the velocity on the boundary to obtain a unique choice of $R$. 

We test this method by applying it to several model systems. In Section \ref{sec:motex} we use a simple example of tethered points to illustrate how the ill-posedness of the free space Stokes equations can lead to nonphysical behavior of the regularized Stokeslet solution. In both this example and the model of a cell migrating through an ECM in Section \ref{sec:cellmot}, we show that subtracting the mean force can create translations in the structure configurations, which for our applications are problematic because we seek measurements on the cell displacement. Because of this, we find our solution of solving the Stokes equations with a zero mean flow on the boundary to give the most physically relevant results for our applications. Finally, we apply this technique to show how force imbalances on the permeable cell cortex drive pressure relief (independent of bleb formation) in blebbing cells. 

We emphasize that this technique is not a solution to generally address Stokes' paradox. The problem remains that no flow is truly 2D, and so representing 3D flows in two dimensions introduces modeling error. However, there are ways to address models with nonzero net forcing so that insight can be gained from 2D models without having to take on the computational complexity of 3D. Here we describe a method to address Stokes' paradox and show that for our applications, the method gives solutions free of artificial translations. Additionally, the method is straightforward to implement and does not involve solving linear systems.

Future work involves extending the work of Section \ref{sec:cellmot} to use our approach to investigate different mechanisms of cell migration (rear contraction in addition to frontal protrusion). We plan to examine the effectiveness of each mechanism for different values of ECM density, ECM stiffness, and cortical tension. 

\section*{Acknowledgements}
The authors thank Ricardo Cortez, Saverio Spagnolie, Charles Peskin, Aleksander Donev, and Yoichiro Mori for helpful discussions. We thank Alex Mogilner for discussions regarding the cell migration model. OM is supported by the NSF Graduate Research Fellowship DGE-1342536 and the Henry MacCracken fellowship. This work was also supported by a grant from the Simons Foundation (\#429808, Wanda Strychalski).
\section*{References}
\bibliographystyle{model5-names}
\bibliography{StokesBib}

\end{document}